\documentclass[12pt]{article}
\usepackage{a4wide}
\usepackage{amsmath}
\usepackage{times}
\numberwithin{equation}{section}
\renewcommand{\thepage}{\bf\arabic{page}}

\catcode`\@=11
\newif\if@fewtab\@fewtabtrue

{\count255=\time\divide\count255 by 60
\xdef\hourmin{\number\count255}
\multiply\count255 by-60\advance\count255 by\time
\xdef\hourmin{\hourmin:\ifnum\count255<10 0\fi\the\count255}}
\def\ps@draft{\let\@mkboth\@gobbletwo
    \def\@oddhead{}
    \def\@oddfoot
       {\hbox to 7 cm{\texttt{\scriptsize{Draft\ version:\ \draftdate}}
       \hfil}\hskip -7cm\hfil\rm\thepage \hfil}
    \def\@evenhead{}\let\@evenfoot\@oddfoot}

\def\mylabel#1{\label{#1}\ifnum\draftcontrol=1\hskip0.5cm
\hbox to2cm{\texttt{\scriptsize{#1}}}\hskip-2cm\else{}\fi}

\def\mybibitem[#1]#2{\bibitem[#1]{#2}\ifnum\draftcontrol=1
\hskip -4cm \hbox to 2cm
{\hfil\texttt{\scriptsize{#2}}}\hskip 2cm\else{}\fi}

\def\draftdate{\number\day/\number\month/\number\year\ \ \ \hourmin }

\global\def\draftcontrol{0}
\catcode`\@=12
\def\eqref#1{(\ref{#1})}

\def\CLAIM#1#2#3{\ifnum\draftcontrol=1\par\hskip -3cm\hbox to
2cm{\hfil\texttt{\scriptsize{#2}}}\hskip 1cm\else{}\fi
\begin{#1}\hskip -5.6pt{\bf .}\hskip 5.6pt\label{#2}#3\end{#1}}
\def\SECT#1#2{
\ifnum\draftcontrol=1\hskip -3cm\hbox to
2cm{\hfil\texttt{\scriptsize{#1}}}\hskip 1cm\else{}\fi
\section{#2}\label{#1}}
\def\PROOF{\medskip\noindent{\bf Proof.\ }}
\def\REMARK#1{\medskip\noindent{\bf Remark.\ }#1\medskip}
\def\LIKEREMARK#1{\medskip\noindent{\bf #1.\ }}
\def\QED{\hfill\vbox{\hrule height 0.4pt
	\hbox{\vrule width 0.4pt height 1.6ex \kern 1.6ex
		\vrule width 0.4pt}
	\hrule height 0.4pt}}
\let\truett=\tt
\def\tt#1{\ifnum\draftcontrol=1\hskip -5truecm{\truett\#\#\#\#\#\#\#\#\#\#\#\#\#\#{#1}}\else{}\fi}
\renewcommand{\AA}{{\cal A}}
\newcommand{\BB}{{\cal B}}
\newcommand{\CC}{{\cal C}}

\newcommand{\EE}{{\cal E}}
\newcommand{\FF}{{\cal F}}
\newcommand{\GG}{{\cal G}}
\newcommand{\HH}{{\cal H}}

\newcommand{\LL}{{\cal L}}
\newcommand{\MM}{{\cal M}}
\newcommand{\NN}{{\cal N}}

\newcommand{\PP}{{\cal P}}

\newcommand{\UU}{{\cal U}}
\newcommand{\VV}{{\cal V}}

\newcommand{\XX}{{\cal X}}

\renewcommand{\natural}{{\bf N}}
\newcommand{\integer}{{\bf Z}}

\newcommand{\real}{{\bf R}}
\newcommand{\complex}{{\bf C}}
\renewcommand\Re{\mathrm{Re\,}}
\renewcommand\Im{\mathrm{Im\,}}
\def\HALF{\frac{1}{2}}

\def\FOUR{\frac{1}{4}}

\let\phi=\varphi
\let\epsilon=\varepsilon

\let\tilde=\widetilde
\def\to{\rightarrow}

\def\lpi(#1){{\mathrm{L}^{\!#1}_{\phantom{0}}}}
\def\ltwoloc(#1){{\mathrm{L}^{\!2}_{#1}}}
\def\ltwo{{\mathrm{L}^{\!2}_{\phantom{0}}}}
\def\linf{{\mathrm{L}^{\!\infty}_{\phantom{0}}}}
\def\ii{\mathrm{i}}
\def\prob{{\bf P}}
\def\expec{{\bf E}}
\def\htop{h_{\mathrm{top}}}
\def\dup{d_{\mathrm{up}}}
\def\cbinf{\CC_{\mathrm{b}}^\infty}
\def\hul(#1){{\mathrm{H}^{#1}_{\mathrm{ul}}}}
\def\hl(#1)(#2){{\mathrm{H}^{#1}_{#2}}}
\def\ho(#1){{\mathrm{H}^{#1}_{\phantom{0}}}}
\def\d(#1){\,\mathrm{d}#1}

\def\wlim{\mathop{\mathrm{w\!-\!lim}}}
\def\card{\mathop{\mathrm{card}}}
\def\norm(#1){\vert\kern-0.1em\vert\kern-0.1em\vert{#1}\vert\kern-0.1em\vert\kern-0.1em\vert}
\def\bdelta{{\boldsymbol\delta }}
\def\xin{\xi^{(n)}_y}
\def\xino{\xi^{(n)}_0}

\begin{document}
\markboth{\vspace{1cm}\bf{Stochastic Ginzburg--Landau Equations}}
{\bf{Stochastic Ginzburg--Landau Equations}}

\title{\bf{Space-Time Invariant Measures, Entropy, and Dimension for
Stochastic Ginzburg--Landau Equations}} 

\author{\vspace{5mm} Jacques Rougemont\\
\vspace{-4pt}\small{Department of Mathematics, Heriot--Watt
University, Edinburgh EH14 4AS, United Kingdom.}}

\date{}

\maketitle

\vspace{5mm}
\begin{center}{\small\texttt{J.Rougemont@ma.hw.ac.uk}\\}\end{center}
\vspace{5mm}
\vfill

\thispagestyle{empty}
\begin{abstract}\noindent We consider a randomly forced
Ginzburg--Landau equation on an unbounded domain. The forcing is smooth
and homogeneous in space and white noise in time. We prove existence
and smoothness of solutions, existence of an invariant measure for the
corresponding Markov process and we define the spatial densities of
topological entropy, of measure-theoretic entropy, and of upper
box-counting dimension. We prove inequalities relating these
different quantities. The proof of existence
of an invariant measure uses the compact embedding of some space of
uniformly smooth functions into the space of locally square-integrable
functions and a priori bounds on the semi-flow in these spaces. 
The bounds on
the entropy follow from spatially localised estimates on the rate of
divergence of nearby orbits and on the smoothing effect of the evolution.
\end{abstract} 
\break

\setcounter{page}{1}

\SECT{intro}{Introduction}The use of dynamical system techniques and
ideas in the study of extended partial differential equations has proved
extremely fruitful in the past, see for example P. Collet's talk at
ICM'98, \cite{colletICM} (where he also emphasises the limitations of
such an approach).
However, until now only results using topological or
geometric properties of the dynamics have been used (like invariant
manifolds, bifurcation theory, topological entropy, Hausdorff
dimension). That is to say, extended dynamical systems are usually
regarded as topological dynamical systems. 
In contrast, most of the very deep results in
finite dimensional dynamical systems use measure-theoretic ideas,
namely ergodic theory (as advocated for instance in the review by
L.-S. Young at the ICMP in 1997, \cite{young}). 

One of the favourite models of infinite dimensional
dynamical systems studied recently is the Ginzburg--Landau equation. It
appears as a generic normal form describing the amplitude of periodic
bifurcated solutions (see \cite{colletICM}) and it is also
believed to be a good example of spatio-temporally chaotic dynamics
\cite{lever}. It is known that its attractor is
infinite-dimensional and has positive $\epsilon $--entropy (see
\cite{coleck1,coleck2,coleck3,rougemont}).

Here we propose to use random perturbations to obtain,
by probabilistic techniques, the existence of invariant measures for
the corresponding random dynamical system. The existence result is
based on the observation by J. Ginibre and G. Velo that the
Ginzburg--Landau equation has global solutions 
both in uniformly local Sobolev spaces and in local $\ltwo$
space. Their proofs go through to the stochastic case without much
effort, if we assume the noise to be smooth in space. Since
uniformly local Sobolev spaces of sufficiently high order are
compactly embedded into local $\ltwo$ space,
we get the Feller property of the semi-group and the tightness of the
Ces\`aro means therefore existence of an invariant measure 
(by standard arguments for stochastic differential
equations, see \cite{daprato}). These measures are also
translation invariant, because the noise, the deterministic
part of the equation and the spaces used are all translation invariant. We
refer to the property of being invariant under the time evolution as
``stationarity'' and the invariance under space translations as
``homogeneity'' of the measure, following Vishik and Fursikov
\cite{vishik}.

In a second part of the paper 
we define the topological entropy and the measure-theoretic
entropy, or rather their spatial densities, since both quantities are
extensive (this has been discovered in this context by Collet and
Eckmann in \cite{coleck1}, see \textit{e.g.\ }\cite{ruelle} for
earlier similar ideas). Usual inequalities from ergodic theory can be
proved in this case and the Collet-Eckmann bound on the topological
entropy is also valid (see \cite{coleck2}).

The paper is organised as follows: In Section~\ref{model}, we set the
model and the functional analysis background needed for the remainder
of the paper. The main results of the paper are summarised in
Section~\ref{informal}.  
In Section~\ref{smoothsol} we obtain uniform bounds on the
solutions in Sobolev spaces, these bounds being then used in
Section~\ref{invariant} to prove the existence of invariant measures. 
Section~\ref{entropy} is devoted to the results on existence and
properties of the (measure-theoretic and
topological) entropies. Various technical proofs have been relegated to
Sections~\ref{iteration}--\ref{gagliardo-proof}. 

We finish this introduction by commenting on the fact that many new
results on invariant measures for nonlinear PDEs have recently
appeared. We mention for instance
\cite{bkl,flandoli,hairer,armen,kuksin,mattingly,sinai}. 
To the best of our
knowledge the present work is the first where the model considered enjoys:
infinite volume (hence continuous spectrum without gap), genuinely
nonlinear interaction, homogeneous noise (hence infinite supply of
energy at each time), non-trivial deterministic dynamics 
({\it e.g.\ }the attractor of the deterministic Ginzburg--Landau is
infinite dimensional). However, we are
still unable to prove uniqueness of the invariant measure ({\it i.e.\ }an 
ergodicity result, as for example in
\cite{armen,funaki1,flandoli,mattingly,bkl,daprato2}).

\LIKEREMARK{Acknowledgements}This work was supported by the Fonds
National Suisse. I am grateful to Sergei Kuksin, Armen Shirikyan, and
Martin Hairer for their comments and suggestions.

\SECT{model}{Model and Definitions}
We consider equations of the form 
\begin{equation}
\begin{split}
&\d(u)\,=\,\Bigl((1+\ii\alpha )\Delta u+u
-(1+\ii\beta )|u|^{2q}u\Bigr)\d(t)+\xi 
\d(w)(t)~,\\
&u(x,t)\,\in\,\complex~,\qquad x\,\in\,\real^d~,
\qquad t\,\ge\,0~,\qquad\alpha ,\beta \,\in\,\real~,
\end{split}
\mylabel{theeq}
\end{equation}
where $w(t)$ is a Wiener process and Eq.\eqref{theeq} is
understood as an It\^ o stochastic differential in $t$. For a while we
simply assume $\xi(\cdot,t)\in\cbinf(\real^d)$ uniformly in $t$ 
and it is adapted to the Wiener process in $t$. A specific example will be
considered in Section~\ref{invariant}. We also assume
$u(\cdot,0)=u_0\in\cbinf(\real^d)$. We make the following
\CLAIM{Hypothesis}{Hn-fct}{We assume $d\le 2$, $q>\HALF$, and
\begin{equation}
-(1+\alpha \beta)\,<\,|\alpha-\beta|\frac{\sqrt{2q+1}}{q}~,\qquad
|\beta |\,\le\,\frac{\sqrt{2q+1}}{q}~.
\mylabel{conditions}
\end{equation}}

\REMARK{The second inequality in \eqref{conditions} implies the first
one. We wrote the first one because it appears in this form in the
proof of Proposition~\ref{Phone-bounds}, while the second
condition appears in the proof of Lemma~\ref{Lextension},
see Ginibre and Velo \cite{ginibre1,ginibre2} for the most general
results available in this direction.}

We next introduce the function spaces used in this paper. Let
\begin{equation}\begin{split}
\phi_{\delta,y}(x)\,=\,\exp\left(-\sqrt{1+\delta ^2|x-y|^2}\right)~.
\end{split}
\mylabel{cutoff}
\end{equation}
The main feature of this function is that it belongs to $\lpi(p)$ for
all $p$ and 
\begin{equation}\begin{split}
\left \|\frac{\nabla^n\phi_{\delta,y}}{\phi_{\delta,y}}\right\|_\infty
\,=\,A_n\delta^n\,<\,\infty
\end{split}
\mylabel{cutoff-bound}
\end{equation}
for all $y\in\real^d$ and $n\in\natural$. This function is used as
weight on Sobolev and Lebesgue spaces:
\CLAIM{Definition}{Dsobul}{The local Lebesgue space 
$\ltwoloc(\delta ,y)$ 
is defined as the completion of $\cbinf$ 
(bounded smooth functions) in the norm induced by the scalar product 
$$
\bigl (f,g\bigr )_{\delta ,y}\,=\,\int\phi _{\delta,y}(x)
\overline{f(x)}g(x)\d(x)~.
$$
The local Sobolev spaces $\hl(m)(\delta ,y)$ are defined as
$$
\hl(m)(\delta ,y)\,=\,\bigl\{f~:~\nabla^kf\,\in\,
\ltwoloc(\delta ,y)~,~k=0,\dots,m\bigr\}~.
$$
The uniformly local Sobolev spaces $\hul(m)$ are
defined as the completion of $\cbinf$ in the norm
$$
\|f\|^2_\hul(m)\,=\,\sum_{k=0}^m\sup_{y\in\real^d}
\bigl (\nabla^kf,\nabla^kf\bigr )_{\delta,y}~.
$$}

Remark that $\hul(m)$ is actually independent of $\delta >0$, since
the following inclusion holds:
\begin{equation}\begin{split}
\hul(m)\,=\,\bigcap_{y\in\real^d}\hl(m)(\delta ,y)
\,\subset\,\bigcap_{\delta>0}\hl(m)(\delta ,y)~.
\end{split}
\mylabel{UB40}
\end{equation}
Usual Sobolev embeddings hold \cite{adams}, for example if $m>d/2$, the
inequality
\begin{equation}\begin{split}
\|f\|^2_\infty\,\le\,C\delta \|f\|^2_\hul(m)
\end{split}
\mylabel{sobolev-ineq}
\end{equation}
implies the continuous embedding $\hul(m)\hookrightarrow\linf$. 
Moreover, by the Rellich--Kondrachov Theorem \cite{adams},
\begin{equation}\begin{split}
\hl(m+k)(\delta ,y)\hookrightarrow\hl(m)(\delta ',y)
\end{split}
\mylabel{rellich-kondrachov}
\end{equation}
is compact if $k>d/2$ and $0<\delta <\delta '$ (see
Section~\ref{compact}).

\LIKEREMARK{Notations}Throughout the paper, $\overline{z}$ denotes the
complex conjugate of $z$, $f_t(x)\equiv f(x,t)$ hence $\|f_t\|_\XX$ is
the norm of $f(x,t)$ in the space $\XX(\!\d(x))$ ({\it e.g.\
}$\XX=\ltwo$). Norms in Lebesgue spaces $\lpi(p)$ are denoted
$\|\cdot\|_p$ ($\|\cdot\|$ is usually the norm on $\ltwoloc(\delta ,y)$
for the current choice of $\delta $ and $y$). 
Expectations and probabilities with respect to the Wiener measure are
denoted $\expec$  and $\prob$. We denote the integer part of
the positive real $x$ by $[x]\equiv\max\{n\in\natural:n\le
x\}$. Symbols $C,C_1,C_2,\dots,c,c_1,c_2,\dots $ usually denote
generic numerical constant. The product $f\star g$ means the
convolution of the functions $f$ and $g$. The cube of side $L$ and
centre $0$ in $\real^d$ is $Q_L=[-\HALF L,\HALF L]^d$.

\SECT{informal}{Summary of Results}

In this section, we describe in a rather informal way the main results of
this paper. The first result (Section~\ref{smoothsol}) is the
following theorem of existence of smooth bounded solutions to
Eq.\eqref{theeq}:
\noindent{\bf Theorem A.~\,}{\it If Hypothesis~\ref{Hn-fct} holds, 
then Eq.\eqref{theeq} with initial data $u_0\in\cbinf$ 
has a unique solution $u(x,t)=u_t(x)$. For all real $p\ge 1$ and integer
$m$, there is a $B_{p,m}<\infty $ such that for all $t>0$:
$$
\expec\|u_t\|^p_\hul(m)\,\le\,B_{p,m}~.
$$}\smallskip

The proof relies on well-known estimates
\cite{ginibre1,mielke,collet-therm} using the dissipative nature
of the nonlinear 
term in Eq.\eqref{theeq} for the deterministic part and on It\^ o's
Lemma to treat the stochastic term. Actually, in the evolution 
equation for
$\expec\|u_t\|^p_{\hul(m)}$, It\^ o's Lemma only generates terms which
are dominated by the nonlinear dissipative term and this implies that
the techniques which were developed for the deterministic equation are
applicable.

By Lemma~\ref{Lextension} we can extend the existence result to the
space $\ltwoloc(\delta ,y)$, hence we can define the Markovian Feller
semi-group $\PP_t$ acting on  $\CC_{\mathrm{b}}
(\ltwoloc(\delta,y),\complex)$ 
for any specific choice of $y$, for example $y=0$:
\begin{equation*}
\begin{split}
\bigl(\PP_t f\bigr)(u)\,=\,\int_{\ltwoloc(\delta ,0)}f(\eta )
\prob\bigl(u_t\,\in\d(\eta )\bigr)~.
\end{split}
\end{equation*}
An invariant measure for Eq.\eqref{theeq} is a fixed point of the dual
semi-group $\PP_t^*$. We next assume that $\xi$ is an homogeneous
process adapted to the Brownian motion. Since $\hul(m)$ is compactly
embedded into $\ltwoloc(\delta ,y)$ for $m>d/2$ (see
\eqref{rellich-kondrachov}) the following theorem (see
Section~\ref{invariant}) is an immediate consequence of Theorem~A by
the Prokhorov and Krylov--Bogolyubov Theorems (see
\cite{arnold,vishik,daprato2}):

\noindent{\bf Theorem B.~\,}{\it There exists at least one invariant
measure $\mu$ for Eq.\eqref{theeq}. This measure is homogeneous in $x$
and its support is contained in $\bigcap_{m\ge 0}\hul(m)$.}\smallskip

Finally, in Section~\ref{entropy}, we define the random attractor (see
\cite{crauel})
\begin{equation*}
\begin{split}
\AA_\omega\,=&\,\overline{\bigcup_{R>0}\AA(\omega ,R)}^{\hul(m)}~,\\
\AA(\omega ,R)\,=&\,
\bigcap_{T>0}\overline{\bigcup_{t>T}
\Phi^t_{\theta^{-t}\omega}(B_R)}^{\hul(m)}~.
\end{split}
\end{equation*}
Here and below $\Phi_\omega ^t$ is the semi-group generated by
Eq.\eqref{theeq}, $\theta^t$ is the time-shift of the noise, and $T_x$
the group of spatial translations. Moreover $B_R\subset\hul(m)$ is the
ball of radius $R$ and centre $0$.
We introduce the following dynamical observables (see
\cite{katok,liu}):
\begin{equation}
\begin{split}
\htop\,=&\,
\lim_{\epsilon \to 0}\lim_{L\to\infty}\frac{1}{L^d}
\lim_{n\to\infty }\frac{1}{n\tau}
\int\log N_{\omega ,n,\tau,Q_L,\epsilon }\,\prob(\!\d(\omega) )~,\\
h_\mu\,=&\,\lim_{\epsilon \to 0}
\lim_{L\to\infty}\frac{1}{L^d}
\lim_{n\to\infty}\frac{1}{n\tau}\int
H_\mu\left(\bigvee_{x\in\integer^d\cap Q_L}~
\bigvee_{k=0}^{n-1}\Phi_\omega^{-k\tau}T_{-x}
(\Sigma_{\theta^{k\tau}T_x\omega,\epsilon})
\right)\,\prob(\!\d(\omega ))~,\\
\HH_\epsilon
\,=&\,\lim_{L\to\infty }\int\frac{\log M_{\epsilon,Q_L,\omega}}{L^d}
\,\prob(\!\d(\omega ))~,\\
\dup\,=&\,\limsup_{\epsilon \to 0}
\frac{\HH_\epsilon}{\log\epsilon ^{-1}} ~,
\end{split}
\mylabel{entropies}
\end{equation}
where $N_{\omega ,n,\tau,Q,\epsilon }$ is the cardinality of a minimal
$(n,\epsilon )$--cover of $\AA_\omega|_Q $, 
$\Sigma_{\omega ,\epsilon }$ is a sequence of 
partitions of $\AA_\omega $ in sets of diameter at most $\epsilon$ 
in the metric of $\linf(Q_1)$, 
$M_{\epsilon ,Q,\omega}$ is the
least cardinality of an $\epsilon $--cover of $\AA_\omega|_Q$ 
and $Q_L=[-\HALF L,\HALF L]^d$ (see
Section~\ref{entropy} for detailed definitions).

The quantities in Eq.\eqref{entropies} are called respectively the
topological entropy, the metric or mea\-sure-theoretic entropy
\cite{katok}, the Kolmogorov--Tikhomirov $\epsilon $--entropy
\cite{kolmogorov} and the upper (box-counting) dimension.
It is important to note that the above numbers are all spatial
densities (limit as $L\to\infty $ of quantities divided by
$L^d$) although the limits are not taken in the most natural order. 
They are thus spatially localised versions of the usual
entropies and dimensions.

We then prove the following estimates:

\noindent{\bf Theorem C.~\,}{\it There is a $\gamma<\infty $ such that
$h_\mu\le \htop\le\gamma \dup<\infty$.}\smallskip

The proof that all the various limits in Eq.\eqref{entropies} 
exist relies on standard subadditive bounds \cite{katok}. 
The upper bound on $\dup$ follows from
spatially localised estimates of the rate of divergence of nearby
orbits (Lemma~\ref{Lexpdiv}) as well as the smoothing action of the
evolution (see Section~\ref{iteration}, in particular
Lemma~\ref{Lksplit}). It is similar to the proof of the
deterministic case \cite{coleck2}.

\SECT{smoothsol}{Bounded Smooth Solutions}
Our first result in this paper is the
existence (and uniqueness by Lemma~\ref{Lcontraction}) 
of smooth bounded solutions to Eq.\eqref{theeq}.
\CLAIM{Theorem}{Tgoodbound}{If Hypothesis~\ref{Hn-fct} holds, 
then Eq.\eqref{theeq} with initial data $u_0\in\cbinf$
has a unique solution $u(x,t)=u_t(x)$. For all real $p\ge 1$ and integer
$m$, there is a $B_{p,m}<\infty $ such that for all $t>0$:
$$
\expec\|u_t\|^p_\hul(m)\,\le\,B_{p,m}~.
$$}

\REMARK{This proof is amply simplified by our
assumptions on the regularity of $\xi_t$ in
Eq.\eqref{theeq}. A much more general theory of stochastic PDEs on
$\real^d$ can be found, for example, in Krylov \cite{krylov}. 
Funaki \cite{funaki1,funaki2} has studied a similar equations with
stronger assumptions on the nonlinearity and Eckmann--Hairer
\cite{hairer} have recently proved a
similar result for stochastic forcings with finite energy.}

\PROOF In the first part of the proof, we fix
$y\in\real^d$ and $\delta >0$ such that
$A_1\delta +A_2\delta ^2<1$ (see Eq.\eqref{cutoff-bound}). 
We write $\|\cdot\|$ and $(\cdot,\cdot)$ for the norm and scalar
product in the corresponding space $\ltwoloc(\delta ,y)$.
All bounds will actually turn out to be
uniform in $y$. We stress that scalar products denoted
$(\cdot,\cdot)$ contain the weight $\phi_{\delta ,y}$ (see
Definition~\ref{Dsobul}) hence integration by parts produces commutators 
of the form $\nabla \phi _{\delta,y}/\phi _{\delta ,y}$. From now on,
we also write $\phi $ for $\phi _{\delta ,y}$.

Let $\LL=(1+\ii\alpha)\Delta +1$. 
For $f\in D_m(\Delta)\subset\hl(m)(\delta ,y)$ 
(the domain of the closure in $\hl(m)(\delta ,y)$
of $\Delta$ with core $\cbinf$), the following holds by
Eq.\eqref{cutoff-bound}:
\begin{equation*}\begin{split}
&\Re\bigl (\nabla^mf,\nabla^m\LL f\bigr )\\
&=\,-\bigl(\nabla^{m+1}f,\nabla^{m+1}f\bigr)
+\bigl(\nabla^mf,\nabla^mf\bigr) 
+\Re\bigl (\nabla\phi\phi^{-1}\nabla^{m}f,
(1+\ii\alpha )\nabla^{m+1}f\bigr )\\
&\le\,
-\HALF\bigl(\nabla^{m+1}f,\nabla^{m+1}f\bigr)
+\Bigl (1+\HALF(1+\alpha^2)\Bigr )
\bigl(\nabla^mf,\nabla^mf\bigr )~.
\end{split}\end{equation*}
Namely $\LL-(1+(1+\alpha ^2)/2)$ 
is a dissipative
operator hence by the Lumer--Phillips Theorem \cite{yosida}, $\LL$
generates a strongly continuous quasi-bounded semi-group $\exp(t\LL)$ on
$\hl(m)(\delta ,y)$, with
\begin{equation}\begin{split}
\|e^{t\LL}\|_{\hl(m)(\delta ,y)\to\hl(m)(\delta ,y)}\,\le\,e^{ct}
\end{split}
\mylabel{semi-group}
\end{equation}
for some $c<\infty$. Remark that we may have chosen $\delta$ such that
$c=1+\epsilon $ for arbitrarily small $\epsilon >0$.

We define mild solutions to Eq.\eqref{theeq}
in $\hl(m)(\delta ,y)$ by the Duhamel formula (and the It\^o integral):
\begin{equation}\begin{split}
z_t\,=\,e^{t\LL}z_0
-(1+\ii\beta )\int_0^te^{(t-s)\LL}|z_s|^{2q}z_s\d(s)
+\int_0^te^{(t-s)\LL}\xi _s\d(w_s)~.
\end{split}
\mylabel{true-int-eq}
\end{equation}
We let $P_M:\real^+\to\real^+$ 
be a smooth cutoff function satisfying $P_M(x)=1$
if $x<M$ and $P_M(x)=0$ if $x>M+1$. We introduce this cutoff into the
nonlinear term above, effectively rendering the nonlinearity uniformly
Lipschitz:
\begin{equation}\begin{split}
\tilde z_t\,=\,e^{t\LL}\tilde z_0
-(1+\ii\beta )\int_0^te^{(t-s)\LL}P_M(|\tilde z_s|)
|\tilde z_s|^{2q}\tilde z_s\d(s)
+\int_0^te^{(t-s)\LL}\xi _s\d(w_s)~.
\end{split}
\mylabel{lip-int-eq}
\end{equation}
We next define the random stopping time $\tau(R)$ by
\begin{equation}\begin{split}
\tau(R)\,=\,\min\bigl\{t\,>\,0~:~\|\tilde z_t\|_\infty
\,\ge\,R\bigr\}~.
\end{split}
\mylabel{tautime}
\end{equation}
We fix arbitrarily a positive number $R<M$, and if $\chi_I$ denotes the
characteristic function of the set $I$, 
we consider the following integral equation for $t<\tau(R)$:
\begin{equation}\begin{split}
u_t\,=
\,e^{t\LL}u_0
-(1+\ii\beta )\int_0^te^{(t-s)\LL}P_M(|u_s|)|u_s|^{2q}u_s
\chi_{\{s\le\tau(R)\}}\d(s)
+\int_0^te^{(t-s)\LL}\xi _s\chi_{\{s\le\tau(R)\}}\d(w_s)~.
\end{split}
\mylabel{stop-int-eq}
\end{equation}

The following is a simple consequence of our construction:
\CLAIM{Lemma}{Ltrue-solutions}{There is almost surely a unique
function $u_t$ satisfying Eq.\eqref{stop-int-eq}, this
function is independent of $M>R$ and it also satisfies
Eq.\eqref{true-int-eq} for $t<\tau(R)$.}

\PROOF See \cite{daprato,kuksin} for the properties of the
stochastic convolution and Section~\ref{uniqueness}
for the contraction argument needed to prove uniqueness.\QED

The remaining part of the proof of Theorem \ref{Tgoodbound} 
follows very closely the paper by Mielke
\cite{mielke} which is based on 
\cite{bgo,collet-therm,ginibre1,ginibre2}. 
We first establish uniform bounds in $\ltwoloc(\delta ,y)$.
\CLAIM{Lemma}{Lltwo-bounds}{For all $\delta >0$ and $p\ge 1$,
there are $C_{0,p}(\delta )$ such that the following bound holds
for all $t>0$ and all $y\in\real^d$:
\begin{equation}\begin{split}
\expec\|u_t\|^p_{\ltwoloc(\delta ,y)}\,\le\,C_{0,p}(\delta )~.
\end{split}
\mylabel{ltwo-const}
\end{equation}}

\PROOF We first estimate the square of the 
norm in $\ltwoloc(\delta ,y)=\ltwo$. By It\^o's formula, we have
\begin{equation}\begin{split}
\d(\|u_t\|^2)\,=&\,-2\|\nabla u_t\|^2\d(t)
-2\Re\bigl (\nabla\phi\phi^{-1}u_t,
(1+\ii\alpha )\nabla u_t\bigr )\d(t)+2\|u_t\|^2\d(t)\\
&\,-2\Re\bigl (u_t,(1+\ii\beta )|u_t|^{2q}u_t\bigr )\d(t)
+\|\xi_t\|^2\d(t)+2\Re(u_t,\xi_t)\d(w_t)\\
\le&\,-\|\nabla u_t\|^2\d(t)
+\bigl (2+(1+\alpha ^2)\bigr )\|u_t\|^2\d(t)
-2\bigl (u_t,|u_t|^{2q}u_t\bigr)\d(t)\\
&\,+\|\xi_t\|^2\d(t)+2\Re(u_t,\xi_t)\d(w_t)\\
\le&\,-\|\nabla u_t\|^2\d(t)
+C(\alpha,\|\xi_t\|,q)\d(t)-\|u_t\|^2\d(t)
+2\Re(u_t,\xi_t)\d(w_t)~.
\end{split}
\mylabel{ito1}
\end{equation}
We integrate this last inequality 
over $t$ and take expectations. By standard
arguments the expectation of the It\^ o integral vanishes (recall that we
consider stopped solutions, Eq.\eqref{stop-int-eq}, see
\cite{daprato}) and we obtain
$$
\expec\|u_T\|^2\,\le\,\expec\|u_0\|^2
-\expec\int_0^T\|\nabla u_t\|^2\d(t)
-\expec\int_0^T\bigl (\|u_t\|^2-C\bigr )\d(t)~.
$$
By Gronwall's inequality, this is
$$
\expec\|u_T\|^2\,\le\,\max\left
\{C_{0,2}~,~(\expec\|u_0\|^2-C_{0,2})e^{-T}+C_{0,2}\right \}~.
$$

For higher powers of the $\ltwo$ norm, we use It\^o's formula again:
$$
\frac{1}{p}\d(\|u_t\|^{2p}_2)\,=\,\|u_t\|^{2p-2}\d(\|u_t\|^2)
+2(p-1)\|u_t\|^{2p-4}\bigl (\Re(u_t,\xi_t)\bigr )^2\d(t)~,
$$
hence (after substituting the estimate \eqref{ito1})
$$
\expec\|u_T\|^{2p}_2\,\le\,\expec\|u_0\|^{2p}_2
-\expec\int_0^T\bigl (\|u_t\|^{2p}-C_{0,2p}\bigr )\d(t)~,
$$
which by Gronwall's inequality gives a uniform bound on $\|u_t\|^p$
for $p>2$. For $p\in[1,2)$, we use Jensen's inequality:
$$
\expec\|u_t\|^p\,\le\,\bigl (\expec\|u_t\|^2\bigr)^{p/2}
\,\le\,C_{0,2}^{p/2}\,=\,C_{0,p}~.
$$

If $u_0$ is uniformly bounded 
and because $\|\xi_t\|^p_{\ltwoloc(\delta ,y)}$ is bounded
uniformly in $y$ and $t$, we obtain the uniform bound in the spaces 
$\ltwoloc(\delta ,y)$ for all $y$
\begin{equation}\begin{split}
\sup_{t>0}~\sup_{y\in\real^d}~\expec\|u_t\|_{\ltwoloc(\delta ,y)}^p
\,\le\,C_{0,p}(\delta )~,
\end{split}
\mylabel{l2uniform}
\end{equation}
which proves Lemma~\ref{Lltwo-bounds}.\QED

\CLAIM{Proposition}{Phone-bounds}{For all $\delta >0$ and $p\ge 1$,
there are $C_{1,p}(\delta )$ such that the following bound holds 
for all $t>0$ and all $y\in\real^d$:
$$
\expec\|u_t\|^p_{\hl(1)(\delta ,y)}\,\le\,C_{1,p}(\delta )~.
$$}

\PROOF We first consider the differential
\begin{equation}\begin{split}
\d(\|\nabla u_t\|^2)
\,=&\,-2\|\Delta u_t\|^2\d(t)
-2\Re\bigl (\nabla\phi\phi^{-1}
\nabla u_t,(1+\ii\alpha )\Delta u_t\bigr )\d(t)
+2\|\nabla u_t\|^2\d(t)\\
&\,+2\Re\bigl (\Delta u_t,(1+\ii\beta )|u_t|^{2q}u_t\bigr )\d(t)
+2\Re\bigl (\nabla\phi\phi^{-1}\nabla u_t,
(1+\ii\beta )|u_t|^{2q}u_t\bigr )\d(t)\\
&\,+\|\nabla\xi_t\|^2\d(t)+2\Re(\nabla u_t,\nabla \xi_t)\d(w_t)\\
\le&\,-2\|\Delta u_t\|^2\d(t)+2\|\nabla u_t\|^2\d(t)
+2\Re\bigl (\Delta u_t,(1+\ii\beta )|u_t|^{2q}u_t\bigr )\d(t)\\
&\,+2\bigl (\sqrt{1+\alpha ^2}\|\Delta u_t\|
+\sqrt{1+\beta ^2}\||u_t|^{2q+1}\|\bigr )\|\nabla u_t\|\d(t)\\
&\,+\|\nabla\xi_t\|^2\d(t)
-2\Re(u_t,\phi^{-1}\nabla(\phi\nabla\xi_t))\d(w_t)~,
\end{split}
\mylabel{first-order}
\end{equation}
and we also compute the following differential that will help us to
cancel out some of the terms above:
\begin{equation}\begin{split}
\frac{1}{q+1}\d(\||u_t|^{q+1}\|^2)\,=&\,
2\Re\bigl (|u_t|^{2q}u,(1+\ii\alpha )\Delta u_t\bigr )\d(t)
+2\||u_t|^{q+1}\|^2\d(t)\\
&\,-2\||u_t|^{2q}u_t\|^2\d(t)
+2q\bigl (\Re\bigl (|u_t|^{q-1}u_t,\xi_t\bigr )\bigr )^2\d(t)
+2\Re\bigl (|u_t|^{2q}u_t,\xi_t\bigr )\d(w_t)\\
\le&\,
2\Re\bigl (|u_t|^{2q}u,(1+\ii\alpha )\Delta u_t\bigr )\d(t)
+2\||u_t|^{q+1}\|^2\d(t)+2q\|\xi_t\|^2\||u_t|^q\|^2\d(t)
\\
&\,-2\||u_t|^{2q}u_t\|^2\d(t)
+2\Re\bigl (|u_t|^{2q}u_t,\xi_t\bigr )\d(w_t)~.
\end{split}
\mylabel{to-the-q}
\end{equation}
We take a convex combination of Inequalities~\eqref{to-the-q} and
\eqref{first-order} (here $\lambda \in[0,1]$):
\begin{equation*}\begin{split}
&\HALF\Bigl (\lambda\d(\|\nabla u_t\|^2)
+\frac{(1-\lambda )}{q+1}\d(\||u_t|^{q+1}\|^2)\Bigr )
\\
&\le\,
\bigl (\lambda\|\nabla u_t\|^2
            +(1-\lambda )\||u_t|^{q+1}\|^2\bigr )\d(t)
+\bigl ((1-\lambda )q\|\xi_t\|^2\||u_t|^q\|^2
+\frac{\lambda}{2}\|\nabla\xi_t\|^2\bigr )\d(t)\\
&\phantom{\le\,}
+\lambda 
\bigl (\sqrt{1+\alpha ^2}\|\Delta u_t\|
        +\sqrt{1+\beta ^2}\||u_t|^{2q+1}\|\bigr )
\|\nabla u_t\|\d(t)\\
&\phantom{\le\,}
+\MM\d(t)+\Re\bigl (-\lambda(u_t,\phi^{-1}\nabla(\phi\nabla\xi_t))
+(1-\lambda )b\bigl (|u_t|^{2q}u_t,\xi_t\bigr )\bigr )\d(w_t)~.
\end{split}
\end{equation*}
The term denoted by $\MM$ is treated separately:
\begin{equation*}\begin{split}
\MM\,=&\,-\bigl(\lambda\|\Delta u_t\|^2
+(1-\lambda)\||u_t|^{2q+1}\|^2\bigr)
+\Re\bigl(1-\ii(\lambda\beta -(1-\lambda )\alpha)\bigr)
\bigl(|u_t|^{2q}u_t,\Delta u_t\bigr)\\
\le&\,
-\epsilon\Bigl (\lambda\|\Delta u_t\|^2
+\frac{(1-\lambda )}{q+1}\||u_t|^{2q+1}\|^2\Bigr )
-2(1-\epsilon )\sqrt{\lambda (1-\lambda)}
\bigl|\bigl(|u_t|^{2q}u_t,\Delta u_t\bigr)\bigr|\\
&\,+\Re\bigl(1-\ii(\lambda\beta -(1-\lambda )\alpha)\bigr)
\bigl(|u_t|^{2q}u_t,\Delta u_t\bigr)\\
\equiv&\,
-\epsilon\Bigl (\lambda\|\Delta u_t\|^2
+\frac{(1-\lambda )}{q+1}\||u_t|^{2q+1}\|^2\Bigr )
+\tilde\MM(\epsilon )~. 
\end{split}\end{equation*}
Under Hypothesis~\ref{Hn-fct}, 
there is an $\epsilon >0$ such that $\tilde\MM(\epsilon )$ 
is negative (see {\it e.g.\ }\cite{ginibre1,mielke,bgo}). The proof
goes as follows: we first remark that integration by parts leads to
$$
\bigl(|u_t|^{2q}u_t,\Delta u_t\bigr)\,=\,
-\left(\frac{\nabla\phi}{\phi}|u_t|^{2q}u_t,
\nabla u_t\right)
-(q+1)\int\phi|u_t|^{2q}|\nabla u_t|^2
\Bigl(1+\frac{q}{1+q}\frac{\overline{u_t}^2}{|u_t|^2}
\frac{\nabla u_t^2}{|\nabla u_t|^2}\Bigr)~.
$$
The last bracket above is of the form $1+z$. Its argument can be
estimated as follows:
$|\arg(1+z)|\le\arcsin|z|=\arcsin\frac{q}{1+q}\equiv\theta$. We plug this into $\tilde\MM$:
\begin{equation*}\begin{split}
\tilde\MM(\epsilon )
\le&\,-(q+1)\bigl|\bigl(|u_t|^{2q},|\nabla u_t|^2\bigr)\bigr|
\biggl\{
\bigl(2(1-\epsilon )\sqrt{\lambda (1-\lambda)}+\cos\theta\bigr)
-|\lambda\beta -(1-\lambda )\alpha|\sin\theta\biggr\}\\
&\,+C(\alpha,\beta ,\lambda ,\epsilon )
\||u_t|^{2q+1}\|\|\nabla u_t\|
\end{split}\end{equation*}
The curly bracket above can be made positive by suitably choosing
$\lambda$ and $\epsilon $, namely we take $\lambda=\cos^2\eta$, we
optimise for $\eta$ and we obtain the following condition, 
which is obviously fulfilled for small
$\epsilon >0$ under Hypothesis~\ref{Hn-fct} (remark
that $1/\tan\theta=\sqrt{2q+1}/q$):
$$
-(1+\alpha \beta) -|\beta -\alpha|/\tan\theta+\epsilon
(2-\epsilon )/\sin^2\theta\,\le\,0~,
$$
(see Ginibre and Velo 
\cite{ginibre1} for this argument, Mielke \cite{mielke} has a
slightly different formulation).
We thus obtain
\begin{equation*}\begin{split}
&\HALF\Bigl (\lambda \d(\|\nabla u_t\|^2)
+\frac{(1-\lambda )}{q+1}\d(\||u_t|^{q+1}\|^2)\Bigr )\\
&\le\,
-\epsilon\Bigl (\lambda \|\Delta u_t\|^2
+\frac{(1-\lambda )}{q+1}\||u_t|^{2q+1}\|^2\Bigr )\d(t)\\
&\phantom{\le\,}
+\Bigl (C_1\|\nabla u_t\|^2+C_2\||u_t|^{q+1}\|^2\Bigr )\d(t)
+\Bigl(C_3\|\xi_t\|^2\||u_t|^q\|^2+C_4\|\nabla\xi_t\|^2\Bigr
)\d(t)\\
&\phantom{\le\,}
+\Bigl (C_5\|\Delta u_t\|+C_6\||u_t|^{2q+1}\|+C_7\Bigr)
\|\nabla u_t\|\d(t)\\
&\phantom{\le\,}
+\Re\Bigl(-\lambda(u_t,\phi^{-1}\nabla(\phi\nabla\xi_t))
+(1-\lambda )\bigl (|u_t|^{2q}u_t,\xi_t\bigr )\Bigr)\d(w_t)
\\
&\le\,
C\d(t)
-\frac{\epsilon}{2}\Bigl (\lambda\|\nabla u_t\|^2
+\frac{(1-\lambda )}{q+1}\||u_t|^{q+1}\|^2\Bigr )\d(t)\\
&\phantom{\le\,}
+\Re\bigl (-\lambda(u_t,\phi^{-1}\nabla(\phi\nabla\xi_t))
+(1-\lambda )\bigl (|u_t|^{2q}u_t,\xi_t\bigr )\bigr )\d(w_t)~,
\end{split}\end{equation*}
thanks to the following obvious inequality:
\begin{equation}\begin{split}
-\|\Delta u_t\|^2
\,\le\,-\rho\|\nabla u_t\|^2+C\rho^2\|u_t\|^2
\end{split}
\mylabel{ghee}
\end{equation}
which holds for all $\rho >0$ and for some $C>0$. As before, we
take expectations, integrate over $t$ and we use Gronwall's inequality
to find out the following bound:
\begin{equation}\begin{split}
&\max\left \{\expec\|\nabla u_T\|^2,\expec\||u_T|^{q+1}\|^2\right\}
\\
&\le\,\max\left
\{C_{1,2}~,~
(\expec\|\nabla u_0\|^2
+\expec\||u_0|^{q+1}\|^2-C_{1,2})
e^{-\epsilon T}+C_{1,2}\right
\}~.
\end{split}
\mylabel{nabla+q+1}
\end{equation}
This and Lemma~\ref{Lltwo-bounds} prove
Proposition~\ref{Phone-bounds}.\QED

We next consider solutions $z(x,t)$ to Eq.\eqref{true-int-eq} 
with bounded initial condition. Proposition~\ref{Phone-bounds} on
$u(x,t)$ implies the following:
\CLAIM{Proposition}{Plinf-bound}{For all $p\ge 1$, 
there is a $C_{\infty,p}$ such that for all
$t>0$, the following holds
\begin{equation}\begin{split}
\expec\|z_t\|^p_\infty \,\le\,C_{\infty,p}~.
\end{split}
\mylabel{linf-bound}
\end{equation}}

\PROOF By the bound~\eqref{sobolev-ineq}, 
if $d=1$ then Proposition~\ref{Phone-bounds} implies
the bound~\eqref{linf-bound} for stopped solutions. If $d=2$ 
we need a bound in $\hl(2)(\delta ,y)$. 
This is easily achieved with the help of a
Gagliardo--Nirenberg inequality which we prove in
Section~\ref{gagliardo-proof}:
\CLAIM{Lemma}{Lgagliardo}{Let $f\in\hul(3)(\real^2)$. For all $K>0$ there
are $C(K)$, $\eta$ such that
\begin{equation*}\begin{split}
&\int\phi _{\delta ,y}(x)\left |\Delta\bigl (|f(x)|^{2q}f(x)\bigr)
\overline{\Delta f(x)}\right |\d(x)\\
&\le\,\frac{1}{K}
\int\phi _{\delta ,y}|\nabla ^3f(x)|^2\d(x) 
+C(K)\left (
\sup_y\int\phi_{\delta,y}(x)|f(x)|^{2(q+1)}\d(x)\right )^\eta~.
\end{split}\end{equation*}}

We use Inequality~\eqref{ghee} (with $u_t$ replaced by $\nabla u_t$),
Lemma~\ref{Lgagliardo}, and the 
estimate~\eqref{nabla+q+1} to bound the time
derivative of $\|\Delta u_t\|^2_{\ltwoloc(\delta,x)}$:
\begin{equation*}\begin{split}
\HALF\d(\|\Delta u_t\|^2_{\ltwoloc(\delta,x)})
\,\le&\,-\HALF\|\nabla^3u_t\|^2_{\ltwoloc(\delta,x)}\d(t)
+(1+(1+\alpha ^2))\|\Delta u_t\|^2_{\ltwoloc(\delta,x)}\d(t)
\\
&\,+C(1+\beta ^2)\sup_y\||u_t|^{q+1}\|_{\ltwoloc(\delta ,y)}
^{2\eta}\d(t)+\HALF\|\Delta \xi_t\|^2_{\ltwoloc(\delta,x)}\d(t)
+(\Delta u_t,\Delta\xi_t)_{\delta,x}\d(w_t)\\
\le&\,\left (-\HALF\d(\|\Delta u_t\|^2_{\ltwoloc(\delta,x)})
+C\right )\d(t)
+(u_t,\phi^{-1}\Delta(\phi\Delta\xi_t))_{\delta,x}\d(w_t) ~,
\end{split}\end{equation*}
where $C$ depends on the parameters in Eq.\eqref{theeq} 
(including $\xi_t$)
and on $\|u_t\|^p_{\ltwoloc(\delta,x)}$, $\||u_t|^{q+1}\|^\eta_{\ltwoloc(\delta,x)}$, 
and $\|\nabla u_t\|^2_{\ltwoloc(\delta,x)}$
which satisfy Bounds~\eqref{ltwo-const} and 
\eqref{nabla+q+1} (or rather some extension of it to deal
with the power $\eta$). By the usual Gronwall inequality, this
proves \eqref{linf-bound} for stopped solutions (Eq.\eqref{stop-int-eq}).

We now choose a very large $n_0\gg C_{0,2}+C_{1,2}+C_{2,2}$ and we let 
$$
\EE_n=\bigl\{u~:~\exists\, t\,<\,\tau(2n)\hbox{ s.t. }
\|u_t\|_\infty \,>\,n_0+n\bigr\}~,
$$
where $\tau$ is the stopping time from Eq.\eqref{tautime}.
By Tchebychev's inequality, by Proposition~\ref{Phone-bounds} and
\eqref{sobolev-ineq} we have
$$
\sum_{n=1}^\infty \prob(\EE_n)\,\le\,\sum_{n=1}^\infty 
n^{-2}\sup_{t\le\tau(2n)}\expec(\|u_t\|_\infty^2)\,<\,\infty ~.
$$
By the Borel--Cantelli Lemma, it means that almost surely only finitely
many of the events 
$\EE_n$ happen, and hence $\|u_t\|_\infty $ remains bounded as
the cutoffs $R$ and $M$ in Eq.\eqref{stop-int-eq} are sent to
infinity. Since
$\tau(R)>0$ a.s.\ for all $R>1$ (by the {\it a priori} 
bound of Lemma~\ref{Ltau-bound}) a uniform bound holds in a small
interval of time and this can be iterated indefinitely. 
This implies uniform boundedness of
$z_t$ and a similar argument holds for $\|z_t\|^p_\infty$, $p>1$.\QED

\LIKEREMARK{Proof of Theorem~\ref{Tgoodbound}}By 
Proposition~\ref{Phone-bounds}, it only
remains to show that $\|\nabla ^mu_t\|^p$ is bounded for $m>1$.
Let $p=2$. We assume that it is true for $m-1$ and we consider 
\begin{equation*}\begin{split}
&\HALF\d(\|\nabla^mu_t\|^2)\\
&\le\,
-\HALF\|\nabla^{m+1}u_t\|^2\d(t)
+\bigl (1+\HALF(1+\alpha ^2)\bigr )
\|\nabla ^mu_t\|^2\d(t)\\
&\phantom{\le\,}
-\Re(1+\ii\beta )\bigl(\nabla^mu_t,\nabla^m(|u_t|^{2q}u_t)\bigr )\d(t)
+\HALF\|\nabla^m\xi_t\|^2\d(t)
+\Re(\nabla^m u_t,\nabla^m\xi_t)\d(w_t)\\
&\le\,
-\HALF\|\nabla^{m+1}u_t\|^2\d(t)
+\Bigl (C_1\|\nabla^m u_t\|^2\d(t)
+C_2\|u_t\|^2_{\hl(m-1)(\delta ,y)}+C_3\Bigr)\d(t)\\
&\phantom{\le\,}
+(-1)^m\Re(u_t,\phi^{-1}\nabla^m(\phi\nabla^m\xi_t))\d(w_t)
~,
\end{split}\end{equation*}
Using \eqref{ghee} (with $u_t$ replaced by $\nabla ^{m-1}u_t$),
Proposition~\ref{Plinf-bound}, and the recursion assumption, 
this can be bounded by:
$$
\HALF\d(\|\nabla^mu_t\|^2)\,\le\,
\HALF\bigl (-\|\nabla^mu_t\|^2+C\bigr )\d(t)
+(-1)^m\Re(u_t,\phi^{-1}\nabla^m(\phi\nabla^m\xi_t))\d(w_t)~.
$$
We then take expectations and integrate:
$$
\expec\|\nabla^mu_t\|^2\,\le\,\expec\|\nabla^mu_0\|^2
-\expec\int_0^t\bigl (\|\nabla^m u_s\|^2-C\bigr )\d(s)
~.
$$
The case $p\neq 2$ is similar.\QED

\SECT{invariant}{Invariant Measures}

We now turn to the problem of
the existence of an invariant measure for the process defined by
Eq.\eqref{theeq}. We construct here an explicit example of a smooth
homogeneous random forcing which admits evident generalisations. 
Let $\xi(x)$ be a $\CC^\infty $ almost periodic
function on $\real^d$. Denoting $T_y\xi(x)=\xi(x+y)$, 
the set
$$
G\,=\,\overline{\bigl\{T_y\xi\,:\,y\in\real^d\bigr\}}^{\,\linf}
$$
is
a compact group which can be endowed with the normalised 
Haar measure $h$. We
denote by $(G,\FF_1,h)$ the corresponding probability space
($\FF_1$ is the sigma-algebra of Borel sets) and by
$\xi_y$ the corresponding random variable. Let next $w_\alpha (t)$ be a
standard Brownian motion (vanishing at $0$) on the probability space
$(\CC_0(\real,\real),\FF_2,{\bf W})$ 
($\FF_2$ is the sigma-algebra generated
by the topology of uniform convergence on compact sets and ${\bf W}$
is the Wiener measure). Let $z_\alpha (t)=y(t)-y(0)$ be a 
continuous process on $(G,\FF_1,h)$ adapted to the filtration
generated by the Brownian motion. We
define the stochastic differential $\xi_{y(0)+z_\alpha
(t)}(x)\d(w_\alpha )(t)$ on the probability space 
$(\Omega,\FF,\prob)=(\CC_0(\real,\real)\times G,\FF_2\times\FF_1,
{\bf W}\times h)$. 
Let $\Xi_\omega(x,t)$ denote $\xi_{y(0)+z_\alpha (t)}(x)
w_\alpha (t)$ where
$\Omega\ni\omega=(\alpha,y(0))\in\CC_0(\real,\real)\times G$. By the
nature of Haar measures, $\prob$ is homogeneous in $x$, 
{\it i.e.\ }$T_y^*\prob=\prob$ for all $y\in\real^d$ (see Vishik and
Fursikov \cite{vishik} for a discussion of homogeneous measures).

Let $\Phi^t_\omega $ be the semi-flow generated by
Eq.\eqref{theeq} with noise $\Xi_\omega (x,t)$. 
Using Lemma~\ref{Lextension}, we can define a Markov semi-group
$\PP_t$ acting on $\CC_{\mathrm{b}}(\ltwoloc(\delta ,0),\complex)$ by
\begin{equation}
\begin{split}
\bigl(\PP_t f\bigr)(u)\,=\,\int_{\ltwoloc(\delta ,0)}f(\eta )
\prob\bigl(\Phi_\omega^t(u)\,\in\d(\eta )\bigr)~.
\end{split}
\mylabel{markov}
\end{equation}
$\PP_t$ is a Markovian Feller semi-group (the Feller property follows
from the continuity of $\Phi_\omega ^t$). Its dual $\PP^*_t$ acts on
probability measures over $\ltwoloc(\delta ,0)$ by 
\begin{equation}
\begin{split}
\bigl(\PP^*_t\mu\bigr)(B)\,=\,\int_{\ltwoloc(\delta ,0)}
\prob\bigl(\Phi_\omega^t(u)\,\in\, B\bigr)\mu(\!\d(u))~.
\end{split}
\mylabel{markov-star}
\end{equation}
We call $\mu$ an invariant measure for Eq.\eqref{theeq} if
$\PP^*_t\mu=\mu$ for all $t>0$ (see Arnold \cite{arnold}).

In this section, we prove the following Theorem, which is actually a
simple consequence of the bounds derived in
Section~\ref{smoothsol}.
\CLAIM{Theorem}{Tinvariant}{There exists at least one invariant
measure $\mu$ for Eq.\eqref{theeq}. This measure is homogeneous in $x$
and its support is contained in $\bigcap_{m\ge 0}\hul(m)$.}

\PROOF We consider the family of measures $\{\overline\mu_t\}_{t>0}$,
where
$$
\overline\mu_t\,=\,\frac{1}{t}\int_0^t\PP^*_s\bdelta_0\d(s)~,
$$
$\bdelta_0$ being the unit mass at $0\in\ltwoloc(\delta ,0)$.
By Theorem~\ref{Tgoodbound}, 
this family is tight in $\ltwoloc(\delta ,0)$ for all $\delta >0$. 
Namely for any $\epsilon >0$ there is a compact 
$K_\epsilon\subset\!\subset\ltwoloc(\delta ,0)$
such that $\overline\mu_t(K_\epsilon )>1-\epsilon $. 
For the set $K_\epsilon $ we choose the
ball of radius $R(\epsilon )$ in $\hul(m)$ ($m>d/2$)
for sufficiently large
$R(\epsilon )$ and the compactness follows from \eqref{UB40}
and \eqref{rellich-kondrachov}.
By the Prokhorov Theorem, $\{\overline\mu_t\}_{t>0}$ is weakly
precompact and
thus there is at least one accumulation point $\mu$. By the standard
Krylov--Bogolyubov argument $\mu$ is an invariant measure (see
\cite{arnold,vishik,daprato2} 
for a detailed statement of these procedures).

Let $t_n$, $n=1,2,\dots$ be a sequence such that
$\mu=\wlim_{n\to\infty}\overline\mu_{t_n}$. Let $B_R$ be the ball of
radius $R$ in $\hl(m)(\delta ,y)$ and let 
$f_{y,R}$ be any bounded continuous
function on $\ltwoloc(\delta ,y)$ vanishing 
on $B_R$ (which is a compact set). Since
the topologies of $\ltwoloc(\delta ,0)$ and of $\ltwoloc(\delta ,y)$ are
equivalent, this function is continuous on $\ltwoloc(\delta ,0)$. 
Obviously
$|\int f_{y,R}(\eta )\overline\mu_{t_n}(\!\d(\eta))|<\|f\|_\infty 
\epsilon(R)$ 
for all $n$ and $y$, where $\epsilon (R)\to 0$ as $R\to\infty $. By
weak convergence of $\overline\mu_{t_n}$ to $\mu$ this also holds for
$\mu$ and hence the support of $\mu$ must be contained in
$\bigcap_{y\in\real}\hl(m)(\delta,y)$.

We next prove the homogeneity of $\mu$. Let
$f\in\CC_{\mathrm{b}}(\hul(m),\complex)$ 
and define the translation operator $T_y$ 
by $T_yf(u)=f(T_yu)$. We have
\begin{equation*}
\begin{split}
\int T_yf(\eta )\mu(\!\d(\eta ))
&=\,
\lim_{n\to\infty }\frac{1}{t_n}\int_0^{t_n}
\left(\int f(T_y\eta )\prob(\Phi_\omega ^t(0)\,\in\d(\eta))
\right)\d(t)\\
&=\,\lim_{n\to\infty }\frac{1}{t_n}\int_0^{t_n}
\left(\int f(\eta )\prob(T_y(\Phi_\omega ^t(0))
\,\in\d(\eta ))\right)\d(t)\\
&=\,\lim_{n\to\infty }\frac{1}{t_n}\int_0^{t_n}
\left(\int f(\eta )\prob(\Phi_{T_y\omega} ^t(T_y(0))\,\in\d(\eta))
\right)\d(t)\\
&=\,\lim_{n\to\infty }\frac{1}{t_n}\int_0^{t_n}
\left(\int f(\eta )\prob(\Phi_\omega ^t(0)\,\in\d(\eta))
\right)\d(t)
\,=\,\int f(\eta )\mu(\!\d(\eta ))~,
\end{split}
\end{equation*}
where we have used the homogeneity of $\prob$. Since the above holds
for all $f$, it proves that $\mu$ is homogeneous and the proof of
Theorem~\ref{Tinvariant} is finished.\QED

\REMARK{In the above construction of a tight family of measures, we
could have considered any homogeneous initial measure $\mu_0$ supported by
$\cap_{m\ge 0}\hul(m)$ instead of $\bdelta_0$.}

\SECT{entropy}{Entropy Estimates}In this section, we define and
estimate different notions of entropy for Eq.\eqref{theeq}. We start
with the topological entropy, then the measure-theoretic entropy and
finally the $\epsilon $--entropy. All these quantities are extensive,
hence we actually define their spatial densities. We define the
spatial density of upper box-counting dimension as well.

To do so we first introduce the basic dynamical setup: let
$\Phi_\omega ^t$ ($t>0$) be the solution semi-flow to Eq.\eqref{theeq}
for given noise parameter 
$\omega $ and let $\theta^t$ be the shift semi-flow on
$\Omega$:
$$
\Xi_{\theta^\tau\omega}(x,t)\,=\,
\Xi_\omega(x,t+\tau )-\Xi_\omega(x,t)~.
$$
Let next $S^t$ be
the semi-flow on $\ltwoloc(\delta ,0)\times\Omega$ defined by
\begin{equation*}
\begin{split}
S^t~:~\ltwoloc(\delta ,0)\times\Omega
\,&\to\,\ltwoloc(\delta ,0)\times\Omega\\
(u,\omega)\,&\mapsto\,
\bigl(\Phi _\omega ^t(u),\theta^t(\omega )\bigr)~.
\end{split}
\end{equation*}
We consider the space $\hul(m)$
($m>d$) endowed with the (weaker) topology of uniform
convergence on the compact $Q\subset\!\subset\real^d$. 
By standard embeddings (see \eqref{sobolev-ineq} and 
\eqref{rellich-kondrachov}) bounded
sets of $\hul(m)$ are compact in $\linf(Q)$. 
Following Crauel {\it et al.\ }\cite{crauel}, 
we define the random attractor $\AA_\omega $ as follows: Let $B_R$ be
the ball of radius $R$, centre $0$ in $\hul(m)$ and let 
\begin{equation*}
\begin{split}
\AA_\omega\,=&\,\overline{\bigcup_{R>0}\AA(\omega ,R)}^{\hul(m)}~,\\
\AA(\omega ,R)\,=&\,
\bigcap_{T>0}\overline{\bigcup_{t>T}
\Phi^t_{\theta^{-t}\omega}(B_R)}^{\hul(m)}~.
\end{split}
\end{equation*}
By the estimates of Section~\ref{smoothsol}, $\AA_\omega $ is almost
surely closed
and bounded in $\hul(m)$ hence it is compact in $\linf(Q)$ for any
bounded $Q\subset\real^d$. Moreover the diameter of $\AA_\omega $ in
$\hul(m)$ is less than some $R_\omega$ with $\prob(\omega : R_\omega
<\infty )=1$ and $\expec(\omega \mapsto R_\omega)<\infty $ 
(by Theorem~\ref{Tgoodbound}).
The following equivariance properties hold (we assume
$\theta^tT_x=T_x\theta^t $ for all $(x,t)\in\real^d\times\real^+$):
\begin{equation}
\mylabel{attractor} 
\begin{split}
\Phi_\omega ^t\AA_\omega \,=&\,\AA_{\theta^t\omega }~,\\
T_x\AA_\omega \,=&\,\AA_{T_x\omega }~,
\end{split}
\end{equation}
and it contains the support of any invariant measure for $S^t$. Let
next $\mu$ be an invariant measure in the sense of
Section~\ref{invariant}, namely a stationary measure for the Markov
semi-group \eqref{markov}. We also assume that $\prob$ is an
invariant measure for $\theta^t$. Then $\mu\times\prob$ is an
invariant measure for the dynamical system $(S^t,\XX,\BB)$ where
$\XX=\ltwoloc(\delta ,0)\times\Omega$ and $\BB$ is the associated
sigma-algebra. More precisely, one has
\begin{equation}
\expec\Bigl(\omega \mapsto\bigl(\Phi_\omega^t\bigr)^*\mu\Bigr)\,=\,\mu
\mylabel{equivar}
\end{equation}
(which is only a rephrasing of $\PP_t^*\mu=\mu$, see
Eq.\eqref{markov-star}). We introduce the following definitions:
\CLAIM{Definition}{Dhtop}{Let $\tau>0$, $n\in\natural$, and
$Q\subset\!\subset\real$.
We define a pseudo-metric $d_{\omega ,n,\tau,Q}$ on $\hul(m)$ by
$$
d_{\omega ,n,\tau,Q}(u,v)\,=\,\max_{k=0,\dots,n-1}
\|\Phi_\omega^{k\tau}(u)-\Phi_\omega^{k\tau}(v)\|_{\linf(Q)}~.
$$
Let $N_{\omega ,n,\tau,Q,\epsilon }$ be the cardinality of a minimal
$(n,\epsilon )$--cover of $\AA_\omega|_Q $ (that is
$N_{\omega ,n,\tau,Q,\epsilon }$ is the least number of
open sets whose diameter in the metric $d_{\omega ,n,\tau,Q}$ is at
most $\epsilon $ and whose union contains $\AA_\omega $).}

We define the cube $Q_L=[-\HALF L,\HALF L]^d$. We are now able to
prove the existence of the spatial density of topological entropy
$\htop$ for Eq.\eqref{theeq}:
\CLAIM{Proposition}{Phtop}{For all $\tau>0$ the following limit exists:
\begin{equation}
\htop\,=\,
\lim_{\epsilon \to 0}\lim_{L\to\infty}\frac{1}{L^d}
\lim_{n\to\infty }\frac{1}{n\tau}
\int\log N_{\omega ,n,\tau,Q_L,\epsilon }\,\prob(\!\d(\omega) )~,
\mylabel{htop}
\end{equation}
This limit is independent of $\tau >0$.}
\PROOF The proof is similar to the deterministic case
treated by Collet and Eckmann in 
\cite{coleck2} and is reproduced in Section~\ref{htopproof}.

Let $\UU=\{U_1,\dots,U_k,\dots\}$ be a countable (or finite)
$\mu$--measurable partition of $\AA_\omega$. 
For two partitions $\UU$ and $\VV$, we denote
their refinement $\{U_k\cap V_\ell:U_k\in\UU,
V_\ell\in\VV, \mu(U_k\cap V_\ell)>0\}$ by $\UU\vee\VV$. Moreover 
$\Phi_\omega^{-\tau}(\UU)=\{\Phi_\omega ^{-\tau}(U_k):U_k\in\UU\}$ is
a measurable partition of $\AA_{\theta^{-\tau}\omega }$ 
whenever $\UU$ is a measurable partition of $\AA_\omega $. (Here
$\Phi_\omega ^{-t}$ stands for the inverse of $\Phi_\omega ^t$, namely
$\Phi_\omega ^{-t}(x)$ is the set of all pre-images of $x$.)
\CLAIM{Definition}{Dmuentropy}{Let $H_\mu\big(\UU)$ and $H_\mu(\UU|\VV)$
denote the entropy of a partition and the conditional entropy, both
relative to a given measure $\mu$. They are defined as follows
\begin{equation*}
\begin{split}
H_\mu(\UU)\,&=\,-\sum_{U\in\UU}\mu(U)\log\mu(U)~,\\
H_\mu(\UU|\VV)\,&=\,-\sum_{U\in\UU,V\in\VV}
\mu(U\cap V)\log\left(\frac{\mu(U\cap V)}{\mu(V)}\right)~.
\end{split}
\end{equation*}
We adopt here the convention $0\log 0=0$ therefore
$0<H_\mu(\UU)\le\log\card(\UU)$ (which is possibly infinite for countable
$\UU$).
We also choose an arbitrary sequence $\Sigma_{\omega ,\epsilon }$ 
of partitions of $\AA_\omega $ in sets of diameter at most $\epsilon$ 
in the metric of $\linf(Q_1)$.
}

The second result in this section is the existence of the spatial
density of measure-theoretic entropy $h_\mu$ 
\CLAIM{Proposition}{Phmu}{For all $\tau>0$ the following limit exists:
\begin{equation}
h_\mu\,=\,\lim_{\epsilon \to 0}
\lim_{L\to\infty}\frac{1}{L^d}
\lim_{n\to\infty}\frac{1}{n\tau}\int
H_\mu\left(\bigvee_{x\in\integer^d\cap Q_L}~
\bigvee_{k=0}^{n-1}\Phi_\omega^{-k\tau}T_{-x}
(\Sigma_{\theta^{k\tau}T_x\omega,\epsilon})
\right)\,\prob(\!\d(\omega ))~.
\mylabel{hmu-limit}
\end{equation}
It is independent
of $\tau >0$ and of the particular choice of the
sequence of partitions $\Sigma_{\omega ,\epsilon}$.}
\PROOF Again, the proof is quite standard, see {\it e.g.\
}\cite{katok,liu} or Section~\ref{hmuproof}.

We next introduce the notions of $\epsilon $--entropy $\HH_\epsilon $ 
of Kolmogorov
and Tikhomirov \cite{kolmogorov} and of upper density of 
dimension $\dup$. 
\CLAIM{Definition}{Dkolmogorov}{Let $M_{\epsilon ,Q,\omega}$ be the
least cardinality of an
open cover of $\AA_\omega $ by sets of diameter less than $\epsilon$ in
the metric of $\linf(Q)$ where $Q$ is compact (we call this an
$\epsilon $--cover of $\AA_\omega |_Q$). 
Let $\HH_\epsilon$ be the Kolmogorov--Tikhomirov 
$\epsilon $--entropy defined by
$$
\HH_\epsilon
\,=\,\lim_{L\to\infty }\int\frac{\log M_{\epsilon,Q_L,\omega}}{L^d}
\,\prob(\!\d(\omega ))~,
$$
et let $\dup(\omega )$ be the upper density of dimension of 
$\AA_\omega $:
$$
\dup
\,=\,\limsup_{\epsilon \to 0}\frac{\HH_\epsilon}{\log\epsilon ^{-1}} ~.
$$}

The main results of the section are the following inequalities
involving the different entropies just defined. Corresponding
inequalities in finite dimensional dynamical systems are well-known
\cite{katok}.
\CLAIM{Theorem}{Tentropy-ineq}{There is a $\gamma<\infty $ such that
\begin{equation}
h_\mu\,\le\,\htop
\,\le\,\gamma \dup<\infty ~.
\mylabel{ben-nevis}
\end{equation}}

Before giving the proof of Theorem~\ref{Tentropy-ineq}, we state a
lemma which will prove useful later on.
\CLAIM{Lemma}{Lexpdiv}{There are $C$, $\gamma$ such that for all
(sufficiently large) $L$  and all (sufficiently small) $\epsilon>0$,
if $\|u-v\|_{\linf(Q_L)}\le\epsilon $ then for $t>0$,
$L'=L-C(1+t)\log 1/\epsilon$, one has
$$
\|\Phi_\omega ^t(u)-\Phi_\omega^t(v)\|_{\linf(Q_{L'})}
\,\le\,Ce^{\gamma t}\epsilon
$$
almost surely.}

\PROOF Let $u_t$ and $v_t$ be two solutions to Eq.\eqref{theeq}. By
Lemma~\ref{Lextension}, 
$$
\|u_t-v_t\|_{\ltwoloc(\delta ,0)}\,\le\,
e^{\gamma t}\|u_0-v_0\|_{\ltwoloc(\delta ,0)}
$$
and moreover both $\|u_t\|_\infty$ and $\|v_t\|_\infty $ 
are bounded uniformly in time (see 
Proposition~\ref{Plinf-bound}). Let $K_t(\cdot)$ be the convolution
kernel associated with the semi-group $\exp(t\LL)$ (see
\eqref{semi-group}) and let $r_s=u_s-v_s$. By Duhamel's formula,
\begin{equation*}
\begin{split}
\left|r_t(x)\right|
\,\le&\,\left|K_t\star r_0(x)\right|
+\left|\int_0^tK_{t-s}\star
\bigl(\GG_1(u_s,v_s)r_s
+\GG_2(u_s,v_s)\overline r_s\bigr)(x)\d(s)\right|\\
\le&\,c_1e^{\gamma t}
\Bigl(\epsilon +\sup_{|x-y|^2\le Ct\log1/\epsilon}|r_0(y)|\Bigr)\\
&\,+\sup_{0\le s\le t}
\bigl(\|\GG_1(u_s,v_s)\|_\infty +\|\GG_2(u_s,v_s)\|_\infty \bigr)
\int_0^t\frac{|K_{t-s}|}{\sqrt{\phi _{\delta ,0}}}\star
(\sqrt{\phi _{\delta ,x}}|r_s|)(x)\d(s)\\
\le&\,c_1e^{\gamma t}
\Bigl(\epsilon +\sup_{|x-y|^2\le Ct\log1/\epsilon}|r_0(y)|\Bigr)
+c_2\|\sqrt{\phi _{\delta ,x}}|r_0|\|_2\int_0^te^{\gamma s}
\left\|\frac{|K_{t-s}|}{\sqrt{\phi _{\delta ,0}}}\right\|_2\d(s)\\
\le&\,c_3(1+t)e^{(1+\gamma)t}\Bigl(2\epsilon
+\sup_{|x-y|^2\le Ct\log1/\epsilon}|r_0(y)|
+\sup_{|x-y|\le C\log1/\epsilon}|r_0(y)|
\Bigr)\\
\le&\,4c_3e^{(2+\gamma)t}\epsilon ~,
\end{split}
\end{equation*}
where in the last line we have assumed $|x|\le\HALF
L-C(1+t)\log1/\epsilon$ (hence $|y|\le\HALF L$) and used the assumption
$\sup_{|y|\le L/2}|r_0(y)|\le \epsilon $.

\LIKEREMARK{Proof of Theorem~\ref{Tentropy-ineq}}We split
Theorem~\ref{Tentropy-ineq} into three independent statements, namely
each one of the three inequalities in \eqref{ben-nevis}.
\bigskip

\LIKEREMARK{Proof of $h_\mu\le \htop$}We follow the most
standard proof (originally by Misiurewicz, quoted in
\cite{katok}). We modify the partition
$\Sigma_{\omega,\epsilon }=\{\sigma_1,\dots,\sigma_N\}$ by
``shrinking'' each element, namely by replacing each $\sigma_k$ by a
closed set $U_k$ with $U_k\subset \sigma_k$ and we define 
$U_0=\AA_\omega \backslash\cup_{k=1}^NU_k$. 
We thus obtain a new partition 
$\UU_{\omega,\epsilon } =\{U_0,\dots,U_N\}$ 
and an open cover
$\VV_{\omega,\epsilon }=\{U_1\cup U_0,\dots,U_N\cup U_0\}$. 
We assume that the $U_k$ have been chosen such that
$\int H_\mu(\Sigma_{\omega,\epsilon } |\UU_{\omega,\epsilon } )
\prob(\!\d(\omega ))<1$. Remark that 
$$
\card\Bigl(\bigvee_{x\in\integer^d\cap Q_L}
\bigvee_{j=0}^{n-1}\Phi_\omega ^{-j\tau}T_{-x}
(\UU_{\theta^{-j\tau}T_x\omega,\epsilon  })\Bigr)
\,\le\,2^{nL^d}
\card\Bigl(\bigvee_{x\in\integer^d\cap Q_L}
\bigvee_{j=0}^{n-1}\Phi_\omega ^{-j\tau}T_{-x}
(\VV_{\theta^{-j\tau}T_x\omega,\epsilon})\Bigr)~,
$$
and by Definition~\ref{Dmuentropy}
\begin{equation*}
\begin{split}
&H_\mu\Bigl(\bigvee_{x\in\integer^d\cap Q_L}\bigvee_{j=0}^{n-1}
\Phi_\omega^{-j\tau}T_{-x}
(\UU_{\theta^{-j\tau}T_x\omega,\epsilon})\Bigr)\\
&\le\,\log\card\Bigl(\bigvee_{x\in\integer^d\cap Q_L}\bigvee_{j=0}^{n-1}
\Phi_\omega^{-j\tau}T_{-x}
(\UU_{\theta^{-j\tau}T_x\omega,\epsilon})\Bigr)\\
&\le\,\log\card\Bigl(\bigvee_{x\in\integer^d\cap Q_L}\bigvee_{j=0}^{n-1}
\Phi_\omega^{-j\tau}T_{-x}
(\VV_{\theta^{-j\tau}T_x\omega,\epsilon})\Bigr)
+nL^d\log 2~.
\end{split}
\end{equation*}
Consequently
\begin{equation*}
\begin{split}
&\lim_{L\to\infty }\frac{1}{L^d}\lim_{n\to\infty }\frac{1}{n\tau}
\int H_\mu\Bigl(\bigvee_{x\in\integer^d\cap Q_L}\bigvee_{j=0}^{n-1}
\Phi_\omega^{-j\tau}T_{-x}
(\UU_{\theta^{-j\tau}T_x\omega,\epsilon})\Bigr)
\,\prob(\!\d(\omega ))\\
&\le\,
\lim_{L\to\infty }\frac{1}{L^d}\lim_{n\to\infty }\frac{1}{n\tau}
\int\log\card\Bigl(\bigvee_{x\in\integer^d\cap Q_L}\bigvee_{j=0}^{n-1}
\Phi_\omega^{-j\tau}T_{-x}
(\VV_{\theta^{-j\tau}T_x\omega,\epsilon})\Bigr)\,\prob(\!\d(\omega ))
+\frac{\log C}{\tau}~.
\end{split}
\end{equation*}
Moreover, the difference between the original partition and the new
one is small, namely:
\begin{equation*}
\begin{split}
&\lim_{n\to\infty }\frac{1}{n\tau}
\int H_\mu\bigl(\bigvee_{j=0}^{n-1}\Phi_\omega^{j\tau}
(\Sigma_{\theta^{-j\tau}\omega,\epsilon})\bigr)\,\prob(\!\d(\omega ))\\
&\le\,
\lim_{n\to\infty }\frac{1}{n\tau}
\int H_\mu\bigl(\bigvee_{j=0}^{n-1}
\Phi_\omega^{j\tau}(\UU_{\theta^{-j\tau}\omega,\epsilon})\bigr)
\,\prob(\!\d(\omega ))
+\frac{1}{\tau }
\int H_\mu(\Sigma_{\omega,\epsilon}|
\UU_{\omega,\epsilon})\,\prob(\!\d(\omega ))~.
\end{split}
\end{equation*}
Since all the above holds for arbitrarily large $\tau>0$ we get 
\begin{equation}
h_\mu\,\le\,\lim_{\epsilon\to 0}\lim_{L\to\infty }\frac{1}{L^d} 
\lim_{n\to\infty }\frac{1}{n\tau}\int\log
\card\Bigl(\bigvee_{x\in\integer^d\cap Q_L}\bigvee_{j=0}^{n-1}
\Phi_\omega^{-j\tau}T_{-x}
(\VV_{\theta^{-j\tau}T_x\omega,\epsilon})\Bigr)
\,\prob(\!\d(\omega ))~.
\mylabel{krok}
\end{equation}
Let next $\delta_{\omega,\epsilon }$ be the Lebesgue number of the cover
$\VV_{\omega,\epsilon}$ 
(namely the largest $\delta _{\omega,\epsilon }>0$ such that every
ball of diameter $\delta_{\omega,\epsilon }$ 
is contained in an element of 
$\VV_{\omega,\epsilon}$). Indeed 
$\delta _{\omega,\epsilon }$ is also the Lebesgue number of 
$\bigvee_{j=0}^{n-1}\Phi_\omega^{j\tau}
(\VV_{\theta^{-j\tau}\omega,\epsilon})$
with respect to the metric $d_{\omega ,n,\tau,Q_L}$. Hence
$$
\card\Bigl(\bigvee_{x\in\integer^d\cap Q_L}\bigvee_{j=0}^{n-1}
\Phi_\omega^{-j\tau}T_{-x}
(\VV_{\theta^{-j\tau}T_x\omega,\epsilon})\Bigr)
\,\le\, M_{\delta _{\omega,\epsilon },Q_L,\omega }~,
$$
and this proves that the r.h.s.\ of \eqref{krok} is less than
$\htop$. 
\bigskip

\LIKEREMARK{Proof of $\htop\le\gamma \dup$}The proof
follows \cite{coleck2}. Let $\rho>0$ be such that $\HH_\epsilon
\le\bigl(\dup+\rho\bigr)\log1/\epsilon $ for all $\epsilon
<\epsilon _0$ and then let $L_0=L_0(\epsilon ,\rho )$ be such that all
$L>L_0$ yield
$$
\left|\int\frac{\log M_{\epsilon ,Q_L,\omega}}{L^d}
\,\prob(\!\d(\omega ))
-\HH_\epsilon \right|\,\le\,\rho ~.
$$

Let $L'=L+C(T+1)\log(1/\epsilon )$ and 
$\epsilon'=C^{-1}\exp(-\gamma T)\epsilon$ (see Lemma~\ref{Lexpdiv}). Let 
an $\epsilon '$--cover of $\AA_\omega|_{Q_{L'}}$ (in the sense of
Definition\ref{Dkolmogorov}) be given. 
Then it is also a $(T/\tau,\epsilon)$--cover (in the sense of
Definition~\ref{Dhtop}), hence
$$
N_{\omega ,T/\tau,\tau,Q_L,\epsilon }\,\le\,M_{\epsilon',Q_{L'},\omega }~,
$$
from which follows
\begin{equation*}
\begin{split}
\htop\,&=\,\lim_{\epsilon \to 0}
\lim_{L\to\infty}\frac{1}{L^d}
\lim_{T\to\infty }\frac{1}{T}
\int\log N_{\omega,T/\tau,\tau,Q_L,\epsilon }\,\prob(\!\d(\omega ))\\
\,&=\,
\lim_{\epsilon \to 0}\lim_{L\to\infty}\frac{1}{L^d}\inf_T\frac{1}{T}
\int\log N_{\omega ,T/\tau,\tau,Q_L,\epsilon }\,\prob(\!\d(\omega ))\\
\,&\le\,\lim_{\epsilon \to 0}
\lim_{L\to\infty}\frac{1}{T}\int
\frac{\log M_{\epsilon',Q_{L'},\omega }}{L^d}\,\prob(\!\d(\omega ))\\
\,&\le\,\lim_{\epsilon \to 0}
\lim_{L\to\infty}\frac{1}{T}
\bigl((\dup+\rho)\log1/\epsilon'+\rho\bigr)~.
\end{split}
\end{equation*}
Since $\log1/\epsilon '=\gamma T+\log(C/\epsilon)$, the limit
$T\to\infty $ and $\rho\to 0$ leaves only $\gamma \dup$ on the
r.h.s.\ above.
\bigskip

\LIKEREMARK{Proof of $\dup<\infty $}We want to prove a bound on
$\HH_\epsilon $ of the form $\HH_\epsilon
\le C\log1/\epsilon $ for small $\epsilon >0$. 
To do so we use iteratively the following bound:
\CLAIM{Lemma}{Literative-bound}{There are $A,B_\omega,C>0$ 
such that for all $L>0$ and sufficiently small $\epsilon >0$, one has
almost surely
\begin{equation}
M_{\epsilon ,Q_L,\omega}
\,\le\,M_{2\epsilon,Q_{L+C},\theta^{-1}\omega}
A^{L^d}B^{1/\epsilon^2}_{\theta^{-1}\omega }~.
\mylabel{xx3}
\end{equation}}

The proof of Lemma~\ref{Literative-bound} is postponed to
Section~\ref{iteration}. 
Let $\epsilon >0$, $L>0$. Remember that there is an $R_\omega$ such that 
$M_{R_\omega ,Q_L,\omega}=1$. Let $T_\omega $
be the smallest integer larger than $(\log 2)^{-1}\log(R_\omega
/\epsilon)$. By iterating $T_\omega $ 
times the bound \eqref{xx3}, we obtain
$$
M_{\epsilon ,Q_L,\omega}
\,\le\,\prod_{n=1}^{T_\omega }
A^{(L+(n-1)C)^d}B^{1/\epsilon^2}_{\theta^{-n}\omega}~,
$$
hence
$$
\HH_\epsilon\,=\,\lim_{L\to\infty }
\int\frac{\log M_{\epsilon ,Q_L,\omega}}{L^d}\,\prob(\!\d(\omega ))
\,\le\,\expec(\omega \mapsto T_\omega )\log A\,\le\,C(\log
1/\epsilon+\log\expec(\omega \mapsto R_\omega )) ~,
$$
(by Jensen's inequality) or, by Definition~\ref{Dkolmogorov},
$$
\dup\,=\,\limsup_{\epsilon \to 0}\frac{\HH_\epsilon}
{\log1/\epsilon }
\,\le\,C~,
$$
which is the bound we wanted to prove. With this inequality, the proof
of Theorem~\ref{Tentropy-ineq} is finished.

\QED
\bigskip

\SECT{iteration}{Proof of Lemma~\ref{Literative-bound}}
We give the proof
for the notationally convenient case $d=1$. Let $u$ and $v$ be two
orbits of Eq.\eqref{theeq} with initial conditions $u_0$ and $v_0$ such
that $u_0$ and $v_0$ belong to $\AA_\omega $. The difference $r=u-v$
satisfies almost surely the equation
\begin{equation}
\partial_tr\,=\,\bigl(1+(1+\ii\alpha )\partial_x^2\bigr)r
+\GG_1(u,v)r+\GG_2(u,v)\overline r~,
\mylabel{diff}
\end{equation}
where we have used the notation of Eq.\eqref{factor}.

Let $\chi(x)$ be a smooth and
monotone function satisfying $\chi(x)=1$ if $x\le 1$
and $\chi(x)=0$ if $x\ge2$. We decompose the kernel of $\exp(t\LL)$ into
a low frequency part and a high frequency part:
\begin{align*}
K^{(-)}_t(x)\,&=\,\frac{1}{2\pi}\int_{-\infty}^\infty\!
e^{\ii px+t(1-(1+\ii\alpha )p^2)}\chi(|p/p^*|)\d(p)~,
\\
K^{(+)}_t(x)\,&=\,\frac{1}{2\pi}\int_{-\infty}^\infty\!
e^{\ii px+t(1-(1+\ii\alpha )p^2)}\bigl (1-\chi(|p/p^*|)\bigr)\d(p)~,
\end{align*}
where $p^*>4$ is a sufficiently large real number. We decompose the 
solutions $r_t(x)$ to Eq.\eqref{diff} accordingly:
\begin{align*}
r_t(x)\,&=\,r_t^{(-)}(x)+r_t^{(+)}(x)~,\\
r_t^{(-)}(x)\,&=\,\bigl (K^{(-)}_t\star r_0\bigr )(x)
+\int_0^t\!\bigl(K^{(-)}_{t-s}\star(\GG_1(u_s,v_s)r_s+\GG_2(u_s,v_s)
\overline r_s)\bigr )(x)\d(s)~,\\
r_t^{(+)}(x)\,&=\,\bigl (K^{(+)}_t\star r_0\bigr )(x)+\int_0^t\!
\bigl(K^{(+)}_{t-s}\star(\GG_1(u_s,v_s)r_s+\GG_2(u_s,v_s)\overline r_s)
\bigr)(x)\d(s)~
\end{align*}

The kernels $K^{(-)}_t$ and $K^{(+)}_t$ have some regularity and
decay properties that we next describe: let the Bernstein class
$B_{R,k}$ be the following set of functions:
\begin{equation}
B_{R,k}\,\equiv\,\left \{f\in\linf\,:\,f{\hbox{ extends to an entire
function} }\, ,\,
|f(z)|\,\le\,Re^{k|\Im z|}\right \}~.
\mylabel{Bernstein}
\end{equation}
We have
\CLAIM{Lemma}{Lksplit}{For all $p^*>4$, $t>\HALF$, $f\in\linf$,
$K^{(-)}_t\star f$ is in $B_{R,2p^*}$ with $R\le 2C_0\|f\|_\infty $.
Moreover, for all $n\in\natural$, there is a $C_n>0$ such that 
\begin{align*}
|K^{(-)}_t(x)|\,&\le\,\frac{C_n}{\sqrt{t}}(1+x^2/t)^{-n}~.\\
|K^{(+)}_t(x)|\,&\le\,\frac{C_n}{\sqrt{t}}e^{-(p^*)^2t/2}(1+x^2/t)^{-n}~.
\end{align*}}

The proof of Lemma~\ref{Lksplit} is omitted, see
\cite{coleck2,rougemont}.

Pick a $2\epsilon$--cover of $\AA_\omega |_{Q_{L+C(\epsilon)}}$ 
(which exists a.s.\ by compactness, see Definition~\ref{Dkolmogorov}) 
and let $u$ and
$v$ belong to one of its elements. Then $r_0=u-v$ 
satisfies $|r_0(x)|\le 2\epsilon $ for $|x|\le\HALF(L+C(\epsilon ))$.
Define
$$
\xin(x)\,=\,\frac{1}{(1+(x-y)^2)^{n/2}}~.
$$
Remark that
Lemma~\ref{Lextension} also holds with $\phi _y$ replaced by
$\xin$ ($n\ge 2$). 
Moreover by reproducing the proof of Lemma~\ref{Lexpdiv} using
the bounds from Lemma~\ref{Lksplit} we obtain (for $|x|\le L/2$):
\begin{equation}
\begin{split}
|r_1^{(-)}(x)|\,&\le\,|K^{(-)}_1\star r_0(x)|+
C\int_0^1\|K^{(-)}_{1-s}/\sqrt{\xino}\|_2\|\sqrt{\xin}r_s\|_2\\
\,&\le\, C\epsilon +2C\epsilon \int_0^1\frac{C_n}{\sqrt{1-s}}
e^{\gamma s}\d(s)\\
\,&\le\,A\epsilon ~,
\end{split}
\mylabel{www1}
\end{equation}
where $A$ depends on $n$ but not on $p^*$ and 
\begin{equation}
\begin{split}
|r_1^{(+)}(x)|\,&\le\,|K^{(+)}_1\star r_0(x)|+
C\int_0^1\|K^{(+)}_{1-s}/\sqrt{\xino}\|_2
\|\sqrt{\xin}r_s\|_2\\
\,&\le\,e^{-(p^*)^2/2}\epsilon+
2C\epsilon \int_0^1\frac{C_ne^{-(p^*)^2(1-s)/2}}
{\sqrt{1-s}}e^{\gamma s}\d(s)\cr
\,&\le\,B(p^*)\epsilon ~,
\end{split}
\mylabel{www2}
\end{equation}
where $B(p^*)\to 0$ as $p^*\to\infty $. We choose $p^*$ so large that
$B(p^*)<\HALF$.

We next use a result of Cartwright (see \cite{kolmogorov}, Eq.(191)):
for all $f$ in the Bernstein class
$B_{R,2p^*}$ (see \eqref{Bernstein}), the following identity holds: 
\begin{equation}
f(x)\,=\,\frac{\sin(8p^*x)}{32(p^*)^2}\sum_{n=-\infty }^\infty
(-1)^nf(x_n)\frac{\sin(4p^*(x-x_n))}{(x-x_n)^2}~,
\mylabel{Cartwright}
\end{equation}
where $x_n=\frac{n\pi}{8p^*}$. Let $f,g$ be in $B_{R,2p^*}$. A simple
application of Eq.\eqref{Cartwright} shows that 
$$
\|f-g\|_{\linf(Q_L)}\,\le\,
C\sup_{|n|\le[4p^*L/\pi]+4Cp^*/(\epsilon\pi)}
|f(x_n)-g(x_n)|+\FOUR\epsilon~.
$$
Hence, among all the functions
in $B_{R_\omega ,2p^*}$ that are bounded by $A\epsilon $ in $[-\HALF
L,\HALF L]$ (by \eqref{www1}, $r_1^{(-)}$ is such a function), 
at most $(4A)^{Cp^*L}(4R_\omega /\epsilon )^{Cp^*/\epsilon }$ 
of them are $\epsilon/2$--separated on $Q_L$. By taking a ball of
diameter $\epsilon$ around each of them, and repeating the operation for
each element of the original $2\epsilon $--cover, we get an 
$\epsilon $--cover of $\Phi_\omega^1(\AA_\omega)|_{Q_L}
=\AA_{\theta^1\omega}|_{Q_L}$. 
The number of elements in this cover is at most 
$$
(4A)^{Cp^*L}(4R_\omega /\epsilon )^{Cp^*/\epsilon }
M_{2\epsilon ,Q_{L+C},\omega}~.
$$
The proof of Lemma~\ref{Literative-bound} is complete.\QED

\SECT{htopproof}{Proof of Proposition~\ref{Phtop}}We follow Collet and
Eckmann's proof \cite{coleck2}, which is itself an adaptation of standard
proofs of existence of the topological entropy, see {\it e.g.\
}\cite{katok} and references therein. The proof of
Proposition~\ref{Phtop} is based on the following inequalities: 
\CLAIM{Lemma}{Lsubad}{For all compacts $Q$, $Q'$, all
$m,n\in\natural$ and $\epsilon>\epsilon'>0$ one has
\begin{align}
N_{\omega ,n,\tau,Q,\epsilon}\,\le&\,N_{\omega,n,\tau,Q,\epsilon'}~,
\mylabel{nonincrease}\\
N_{\omega ,n,\tau,Q\cup Q',\epsilon}\,\le&\,N_{\omega,n,\tau,Q,\epsilon}
N_{\omega,n,\tau,Q',\epsilon}~,
\mylabel{Qsubad}\\
N_{\omega ,n+m,\tau,Q,\epsilon}\,\le&\,N_{\omega,n,\tau,Q,\epsilon}
N_{\theta^{n\tau }\omega,m,\tau,Q,\epsilon}~,
\mylabel{tsubad}
\end{align}
Furthermore for any $\tau'<\tau$ the following inequalities hold:
\begin{equation}
N_{\omega,n,\tau',Q_L,\epsilon}\,\le\,
N_{\omega,n,\tau,Q_{f(L)},g(\epsilon)}\,\le\,
N_{\omega,n,\tau',Q_{f(f(L))},g(g(\epsilon))}~,\mylabel{chappo}
\end{equation}
where $f(L)=L+C(\tau +1)\log\epsilon ^{-1}$
and $g(\epsilon )=c\exp(-\gamma\tau)\epsilon$ with $C,c,\gamma$ 
some constants.
}

Lemma~\ref{Lsubad} implies immediately that the limit in
Eq.\eqref{htop} exists: by subadditivity \eqref{tsubad} and by
invariance of $\prob$ under $\theta^t$, we get that
$$
\Lambda_1\,=\,\lim_{n\to\infty }\frac{1}{n\tau }
\int\log N_{ \omega,n,\tau,Q_L,\epsilon}\,\prob(\!\d(\omega ))
$$
exists, it is non-increasing in $\epsilon $ and by further
subadditivity \eqref{Qsubad}
$$
\Lambda _2\,=\,\lim_{L\to\infty }\frac{1}{L^d}\Lambda _1
$$
also exists and is non-increasing in $\epsilon $
(by \eqref{nonincrease}). Hence the limit in
Eq.\eqref{htop} exists. By \eqref{chappo}, it is independent of
$\tau$.\QED

\LIKEREMARK{Proof of Lemma~\ref{Lsubad}}The inequality 
\eqref{nonincrease} is
obvious from the definitions. We prove \eqref{Qsubad} 
by making the observation that 
if $\{A_1,\dots,A_N\}$ is an $(n,\epsilon )$--cover of 
$\AA_\omega |_Q$ and $\{B_1,\dots,B_M\}$ an $(n,\epsilon )$--cover of
$\AA_\omega |_{Q'}$, then $\{A_j\cap B_k:j=1,\dots,N,k=1,\dots,M\}$
is an $(n,\epsilon )$--cover of $\AA_\omega |_{Q\cup Q'}$.

Similarly if
$\{A_1,\dots,A_N\}$ is an $(n,\epsilon )$--cover of
$\AA_\omega |_Q$ and $\{B_1,\dots,B_M\}$ an $(m,\epsilon )$--cover of
$\AA_{\theta^{n\tau}\omega} |_Q$, then $\{A_j\cap \Phi_\omega ^{-n\tau}
B_k:j=1,\dots,N,k=1,\dots,M\}$ 
is an $(m+n,\epsilon )$--cover of $\AA_\omega |_Q$ which
proves \eqref{tsubad}.

The inequality \eqref{chappo} follows
immediately from Lemma~\ref{Lexpdiv}, since if $D$ is
a set of diameter $g(\epsilon)$ in the metric 
$d_{\omega ,n,\tau,Q_{f(L)}}$ then
$D$ is a set of diameter at most $\epsilon$ in the metric 
$d_{\omega ,n,\tau',Q_L}$.\QED

\REMARK{The topology of $\linf(Q)$ is a simplifying choice (as far as
Eq.\eqref{Qsubad} is concerned), but \cite{coleck3} have demonstrated
that other topologies can be used as well.}

\SECT{hmuproof}{Proof of Proposition~\ref{Phmu}}This proof is, like
the proof of Proposition~\ref{Phtop}, based on subadditive
bounds. We use well-known properties of the function $H_\mu(\cdot)$, see
\cite{katok}, Chapter 4.3 (in particular Proposition 4.3.3). We recall
that $x\mapsto-x\log x$ is concave, hence for any partition
$\UU$ and any $t>0$, the following holds:
$$
\int H_\mu\bigl(\Phi_\omega ^{-t}(\UU)\bigr)\,\prob(\!\d(\omega ))
\,\le\,H_\mu\biggl(\int\Phi_\omega ^{-t}(\UU)\,\prob(\!\d(\omega))
\biggr)\,=\,H_\mu(\UU)
$$
where we have used Eq.\eqref{equivar}. We thus have
\begin{equation*}
\begin{split}
&\int
H_\mu\Bigl(
\bigvee_{k=0}^{n+m-1}
\Phi_\omega^{-k\tau}(\Sigma_{\theta^{k\tau}\omega,\epsilon})
\Bigr)\,\prob(\!\d(\omega ))\\
\,&=\,\int H_\mu\Bigl(
\bigvee_{k=0}^{n-1}
\Phi_\omega^{-k\tau}(\Sigma_{\theta^{k\tau}\omega,\epsilon})
\Bigr)\,\prob(\!\d(\omega ))\\
&\phantom{=\,}
+\int H_\mu\Bigl(
\bigvee_{k=n}^{n+m-1}
\Phi_\omega^{-k\tau}(\Sigma_{\theta^{k\tau}\omega,\epsilon})
\Bigl|\Bigr.
\bigvee_{k=0}^{n-1}
\Phi_\omega^{-k\tau}(\Sigma_{\theta^{k\tau}\omega,\epsilon})
\Bigr)\,\prob(\!\d(\omega ))\\
&\le\,\int H_\mu\Bigl(
\bigvee_{k=0}^{n-1}
\Phi_\omega^{-k\tau}(\Sigma_{\theta^{k\tau}\omega,\epsilon})
\Bigr)\,\prob(\!\d(\omega ))
+\int H_\mu\Bigl(
\Phi_{\omega} ^{-n\tau}\bigvee_{k=0}^{m-1}
\Phi_{\theta^{n\tau}\omega}^{-k\tau}
(\Sigma_{\theta^{(k+n)\tau}\omega,\epsilon})
\Bigr)\,\prob(\!\d(\omega ))\\
&\le\,\int H_\mu\Bigl(
\bigvee_{k=0}^{n-1}
\Phi_\omega^{-k\tau}(\Sigma_{\theta^{k\tau}\omega,\epsilon})
\Bigr)\,\prob(\!\d(\omega ))
+\int\int H_\mu\Bigl(
\Phi_{\omega'} ^{-n\tau}\bigvee_{k=0}^{m-1}
\Phi_\omega^{-k\tau}(\Sigma_{\theta^{k\tau}\omega,\epsilon})
\Bigr)\,\prob(\!\d(\omega'))\,\prob(\!\d(\omega))\\
&\le\int H_\mu\Bigl(
\bigvee_{k=0}^{n-1}
\Phi_\omega^{-k\tau}(\Sigma_{\theta^{k\tau}\omega,\epsilon})
\Bigr)\,\prob(\!\d(\omega ))
+\int H_\mu\Bigl(
\bigvee_{k=0}^{m-1}
\Phi_\omega^{-k\tau}(\Sigma_{\theta^{k\tau}\omega,\epsilon})
\Bigr)\,\prob(\!\d(\omega ))~,
\end{split}
\end{equation*}
namely subadditivity in the time variable. We can prove subadditivity in
the space variable in a similar way. 
Thus the first two limits in Eq.\eqref{hmu-limit} exist. These
limits are monotonically increasing as $\epsilon\to 0$, hence the
third limit is well-defined.

We next prove that the limit is independent of the choice of
$\Sigma_{\omega ,\epsilon}$: let $\Sigma_{\omega,\epsilon}$ and
$\tilde \Sigma_{\omega,\epsilon}$ be two different
sequences, we get (by the Rokhlin inequality)
\begin{equation*}
\begin{split}
&\biggl|
\lim_{L\to\infty }\frac{1}{L^d}
\lim_{n\to\infty}\frac{1}{n\tau}H_\mu
\Bigl(
\bigvee_{x\in\atop\integer^d\cap Q_L}
\bigvee_{k=0}^{n-1}\Phi_\omega^{-k\tau}T_{-x}
(\Sigma_{\theta^{k\tau}T_x\omega,\epsilon})
\Bigr)\\
&\phantom{\biggl|}
\,-\lim_{L\to\infty }\frac{1}{L^d}
\lim_{n\to\infty}\frac{1}{n\tau}H_\mu
\Bigl(
\bigvee_{x\in\atop\integer^d\cap Q_L}
\bigvee_{k=0}^{n-1}\Phi_\omega^{-k\tau}T_{-x}
(\tilde \Sigma_{\theta^{k\tau}T_x\omega,\epsilon})
\Bigr)
\biggr|\\
&\le\,H_\mu(\Sigma_{\omega,\epsilon}|\tilde\Sigma_{\omega,\epsilon})
+H_\mu(\tilde \Sigma_{\omega,\epsilon}|\Sigma_{\omega,\epsilon})
\end{split}
\end{equation*}
and the r.h.s.\ above vanishes as $\epsilon \to 0$
since these sequences generate the whole sigma-algebra of $\AA_\omega $
in this limit.

We prove that Eq.\eqref{hmu-limit} is independent of $\tau$ by using
Lemma~\ref{Lexpdiv} and an argument similar to the one used in
Section~\ref{htopproof}.\QED

\SECT{uniqueness}{Uniqueness of Solutions}
In this section, we use the Contraction Mapping Principle to prove
uniqueness of solutions to Eq.\eqref{stop-int-eq} as well as
estimates on the stopping times Eq.\eqref{tautime} (along the lines
of Da Prato and Zabczyk, \cite{daprato} Chapter 7). We define a
Banach space $\BB_{p,K,T}$ 
of complex-valued predictable processes $u$ on the time
interval $[0,T]$ with norm defined by
$$
\norm(u)_p\,=\,\left (\sup_{0\le t\le T}\expec(e^{-Kt}\|u_t\|_\infty
^p)\right )^{1/p}~.
$$
We prove existence and uniqueness of solutions to
Eq.\eqref{lip-int-eq} in $\BB_{p,K,T}$:
\CLAIM{Lemma}{Lcontraction}{Let $T>0$, $p>1$, and $M>1$. 
For sufficiently large
$K$, there exists a unique solution $u_t\in\BB_{p,K,T}$ 
to Eq.\eqref{lip-int-eq} with initial data $u_0$.}

\PROOF We define the map $\FF:\BB_{p,K,T}\to\BB_{p,K,T}$ by (see
Eq.\eqref{lip-int-eq})
$$
\FF\bigl (X\bigr )_t\,=\,e^{t\LL}u_0
+\int_0^te^{(t-s)\LL}P_M(|X_s|)|X_s|^{2q}X_s\d(s)
+\int_0^te^{(t-s)\LL}\xi_s\d(w_s)~.
$$
If $\expec(\|u_0\|_\infty ^p)<\infty $ then obviously $\FF$ maps
$\BB_{p,K,T}$ into itself. We define 
\begin{equation}
\begin{split}
\NN(|x|^2)\,=&\,-(b+\ii\beta) P_M(|x|)|x|^{2q}~,\\
\GG_1(x,y)\,=&\,\HALF\Bigl(\NN(|x|^2)+\NN(|y|^2)
+(|x|^2+|y|^2)\int_0^1\NN\,'\bigl (t|x|^2+(1-t)|y|^2\bigr )\d(t)
\Bigr)~,\\
\GG_2(x,y)\,=&\,xy\int_0^1\NN\,'\bigl (t|x|^2+(1-t)|y|^2\bigr )\d(t)~.
\end{split}
\mylabel{factor}
\end{equation}
If $X$ and $Y$ are arbitrary elements of $\BB_{p,K,T}$, then 
\begin{equation*}\begin{split}
&\norm({\FF(X)-\FF(Y)})^p\\
&=\,\sup_{0\le t\le T}\!\!
\expec\left(\left \|
\int_0^t\!e^{(t-s)(\LL-K)}\left (
\GG_1(X_s,Y_s)e^{-Ks}\bigl (X_s-Y_s\bigr )
+\GG_2(X_s,Y_s)e^{-Ks}\bigl (\overline{X_s}-\overline{Y_s}\bigr
)\right )\!\d(s)
\right \|^p_\infty \right)\\
&\le\,CM^{2qp}(K-c)^{-p}\norm(X-Y)^p~.
\end{split}\end{equation*}
The map $\FF$ is thus a contraction on $\BB_{p,K,T}$ if
$K>C'M^{2q}$. This proves the existence and the uniqueness of
a solution $u_t$ to Eq.\eqref{lip-int-eq}.\QED

To be able to treat solutions to Eq.\eqref{stop-int-eq} as solutions to
Eq.\eqref{true-int-eq} for some time, we use the following bounds on
the stopping times $\tau(R)$ defined by Eq.\eqref{tautime}:
\CLAIM{Lemma}{Ltau-bound}{There is a $C>0$ such that 
the following holds almost surely for all $R>1$: 
$$
\tau(R)\,\ge\,CR^{-2q}\log R~.
$$}

\PROOF This follows immediately from Lemma~\ref{Lcontraction} 
since we can take any $M>R$.\QED

We next show that solutions of Eq.\eqref{theeq} are also uniquely
defined on $\ltwoloc(\delta ,y)$, using that bounded functions form a
dense subset.
\CLAIM{Lemma}{Lextension}{The semi-flow $\Phi_\omega ^t$ extends
almost surely to a bounded continuous semi-flow on
$\ltwoloc(\delta,y)$ for any $\delta >0$ and $y\in\real^d$.}
\PROOF We apply the
non-propagation estimate of Ginibre and Velo \cite{ginibre1}. Let
$u_0$ and $v_0$ be two functions in $\ltwoloc(\delta ,y)$ and denote
the corresponding solutions to Eq.\eqref{theeq} by $u_t$ and
$v_t$. Their difference $u_t-v_t$ satisfies (almost surely) the
following inequality:
\begin{equation*}
\begin{split}
\HALF\partial_t\|\sqrt{\phi _{\delta,y}}(u_t-v_t)\|_2^2
\,\le&\,(1+\HALF\sqrt{1+\alpha ^2})
\|\sqrt{\phi _{\delta,y}}(u_t-v_t)\|^2_2\\
&\,-\Re(1+\ii\beta )\int\phi _{\delta,y}(\overline u_t-\overline v_t)
\bigl(|u_t|^{2q}u_t-|v_t|^{2q}v_t\bigr)~.
\end{split}
\end{equation*}
By \cite{ginibre1} (Proposition 3.1), 
Hypothesis~\ref{Hn-fct} implies that 
the last term above is negative. We thus get an estimate of the
form $\|u_t-v_t\|_{\ltwoloc(\delta,y)}
\le\exp(ct)\|u_0-v_0\|_{\ltwoloc(\delta ,y)}$

This and Lemma~\ref{Lltwo-bounds} 
prove that $\Phi_\omega ^t$ is uniformly bounded and continuous on
$\ltwoloc(\delta ,y)$ for any $\delta >0$ and $y\in\real^d$ if we
define $u_t=\lim_{n\to\infty }u_t^{(n)}$ where $u_0^{(n)}$ is a
Cauchy sequence of bounded functions approaching $u_0$.\QED

\SECT{compact}{Compact Embedding for Local Spaces}In this section, we
give a proof of Relation \eqref{rellich-kondrachov} which is a trivial
adaptation of \cite{adams}, Theorem~6.53, p.174. More precisely we 
prove the embedding \eqref{rellich-kondrachov} to be
Hilbert--Schmidt. 
Let $\{e_n\}_{n\in\natural}$ be a complete orthonormal basis of
$\hl(m+k)(\delta ,y)$. Let $\{Q_n\}_{n\in\natural}$ be a countable
cover of $\real^d$ by balls of radius $1$. Let $x\in Q_n$, 
let $\alpha\le m$ and define the bounded linear operator $D_x^\alpha $ on
$\hl(m+k)(\delta,y)$ by
$$
D_x^\alpha (u)\,=\,\nabla ^\alpha u(x) ~.
$$
Its norm is (by Sobolev embedding) bounded by
$$
\|D_x^\alpha (u)\|^2_{\hl(m+k)(\delta ,y)}
\,\le\,\max_{0\le\alpha\le m}\sup_{x\in Q_n}|\nabla ^\alpha u(x)|^2
\,\le\,\frac{C}{\inf_{x\in Q_n}\phi _{\delta,y}(x)}
\|u\|^2_{\hl(m+k)(\delta ,y)}~.
$$
By Riesz' Lemma, $D_x^\alpha (\cdot)=(v_x^\alpha
,\cdot)_{\hl(m+k)(\delta ,y)}$ for some vector $v_x^\alpha $ and
$$
\sum_{n=1}^\infty |\nabla ^\alpha e_n(x)|^2
\,=\,\sum_{n=1}^\infty|(e_n,v_x^\alpha)_{\hl(m+k)(\delta,y)}|^2
\,=\,\|v_x^\alpha \|^2_{\hl(m+k)(\delta,y)}~.
$$ 
Thus the Hilbert--Schmidt norm of the embedding map is
$$
\sum_{n=1}^\infty \|e_n\|^2_{\hl(m)(\delta',y)}
\,=\, \sum_{\alpha\le m}\int_{\real^d}
\|v_x^\alpha \|^2_{\hl(m+k)(\delta',y)}\phi _{\delta',y}(x)\d(x)
\,\le\,m\sum_{n=1}^\infty\int_{Q_n}
\frac{C\phi_{\delta',y}(x)}{\inf_{z\in Q_n}\phi_{\delta,y}(z)}\d(x)~,
$$
which is finite whenever $\delta '>\delta $.\QED 

\SECT{gagliardo-proof}{Proof of Lemma~\ref{Lgagliardo}}
The proof can be found in \cite{bgo,mielke} 
and is summarised below. We decompose
the plane into countably many sets 
$Q(m,n)$ of unit area and use the bounds 
$\phi_{\delta ,y}(x)\le \exp(-\delta |x-y|)\le
e\phi_{\delta,y}(x)$. For simplicity we assume $\delta =1$ and we
drop it from our notation (if Lemma~\ref{Lgagliardo} is true for 
$\delta=1$ then it is true for all $\delta >0$ by scaling, possibly with
different constants). We simply write $\int_D f$ for $\int_D
f(x)\d(x)$ for $D\subset\real^2$. We have
\begin{equation}\begin{split}
\int_{\real^2}\phi_y|\Delta(|f|^{2q}f)\overline{\Delta f}|
\,\le\,C\sum_{m,n}e^{-|n|}\int_{Q(m,n)}
|\Delta f||f|^{2q-1}\bigl (|f||\Delta f|+|\nabla f|^2\bigr )~,
\end{split}
\mylabel{decomp}
\end{equation}
where $\bigcup_mQ(m,n)\supset\{x\in\real^2:n-\HALF\le |x-y|\le n+\HALF\}$.
We estimate each summand using H\"older and Gagliardo--Nirenberg
inequalities. For any $p,r$ with $p^{-1}+r^{-1}=1$ and in
particular for $r=1+1/q$ and $p=1+q$, we get:
\begin{equation*}\begin{split}
&\int_{Q(m,n)}|\Delta f||f|^{2q-1}\bigl (|f||\Delta f|+|\nabla f|^2\bigr )
\\
&\le\, c_1\|\Delta f\|_{2p}
\bigl (
\|f\|_{2pq/(p-1)}^{2q}\|\Delta f\|_{2p}
+\|f\|_{2pq/(p-1)}^{2q-1}\|\nabla f\|^2_{4pq/(p+q-1)}
\bigr )
\\
&\le\, c_2\|\Delta f\|_{2p}
\Bigl (
\|f\|_{2pq/(p-1)}^{2q}\|\Delta f\|_{2p}
+\|f\|_{2pq/(p-1)}^{2q-1}
\bigl(\|f\|^{1/2}_{2pq/(p-1)}\|\Delta f\|^{1/2}_{2p}\bigr)^2
\Bigr )\\
&=\,c_3\|\Delta f\|_{2p}^2\|f\|_{2qr}^{2q}\\
&\le\,c_4\|\nabla ^3f\|_2^{2(2q+2)/(2q+3)}
\|f\|_{2(q+1)}^{2(q+1/(2q+3))}\\
&\le\,K^{-1}\|\nabla ^3f\|_2^2+c_5K\|f\|_{2(q+1)}^{4q^2+6q+2}~.
\end{split}\end{equation*}
By summing up all contribution to \eqref{decomp} we arrive at
\begin{equation*}\begin{split}
&\int_{\real^2}\phi_y|\Delta(|f|^{2q}f)\overline{\Delta f}|\\
&\le\,CK^{-1}\sum_{m,n}e^{-|n|}\int_{Q(m,n)}|\nabla ^3f|^2
+C'K\sum_{m,n}e^{-|n|}\left(\int_{Q(m,n)}|f|^{2(q+1)}\right
)^\eta\\
&\le\,\tilde CK^{-1}\int_{\real^2}\phi _y|\nabla ^3f|^2
+C''K\sum_n ne^{-|n|}\left(\sup_y\int_{\real^2}\phi _y|f|^{2(q+1)}\right
)^\eta\\
&=\,\tilde CK^{-1}\int_{\real^2}\phi _y|\nabla ^3f|^2
+C'''K\left(\sup_y\int_{\real^2}\phi _y|f|^{2(q+1)}\right)^\eta~,
\end{split}\end{equation*}
which proves Lemma~\ref{Lgagliardo}.\QED

\end{document}